%% file: main.tex
\documentclass[a4paper]{article}
\usepackage{geometry}
\usepackage{cite}
\usepackage{diagbox}
\usepackage{comment}
\usepackage{wrapfig,lipsum,booktabs}
\usepackage{amsmath,amssymb,amsfonts}
\usepackage[ruled,part]{algorithm}
\usepackage{graphicx}
\usepackage{amsthm}
\input{macros}
\usepackage{url}
\usepackage{pgfplots}
\pgfplotsset{compat=1.18}
\usepackage{algorithm}
\usepackage{algpseudocode}
\usepackage{gensymb}
\usepackage[numbered,framed]{matlab-prettifier}

\usepackage{subcaption}
\usepackage{tikz-cd}
\usepackage{listings}
\usepackage{xcolor}

\definecolor{codegreen}{rgb}{0,0.6,0}
\definecolor{codegray}{rgb}{0.5,0.5,0.5}
\definecolor{codepurple}{rgb}{0.58,0,0.82}

\lstdefinestyle{mystyle}{
  commentstyle=\color{codegreen},
  keywordstyle=\color{magenta},
  numberstyle=\tiny\color{codegray},
  stringstyle=\color{codepurple},
  basicstyle=\ttfamily\footnotesize,
  breakatwhitespace=false,
  breaklines=true,
  captionpos=b,
  keepspaces=true,
  numbers=left,
  numbersep=5pt,
  showspaces=false,
  showstringspaces=false,
  showtabs=false,
  frame = single,
  tabsize=2
}
\newtheorem{theorem}{Theorem}
\newtheorem{definition}{Definition}
\usepackage{textcomp}
\usepackage{xcolor}
\def\BibTeX{{\rm B\kern-.05em{\sc i\kern-.025em b}\kern-.08em
    T\kern-.1667em\lower.7ex\hbox{E}\kern-.125emX}}
\begin{document}
\title{Adaptive Methods for Multiobjective Unit Commitment}
\author{Ece Tevruez and Aswin Kannan\footnote{
The first author is with Nebenan.de and is reachable at ecetevruz@gmail.com. The second author is the corresponding author and is reachable at both aswin.kannan@hu-berlin.de and aswin.kannan@iiitb.ac.in. This work was done when the first author was pursuing her Masters at Humboldt Universitaet zu Berlin, Germany under the guidance of the second author. \\
The second author was a faculty member and Junior Research Group Leader at Humboldt Universitaet zu Berlin till December 2024. He is currently a faculty member at the International Institute of Information Technology, Bangalore, India (from January 2025). He still retains an affiliation to his Junior Group in Berlin.}}

\maketitle



\abstract{
This work considers a multiobjective version of the unit commitment problem that deals with finding the optimal generation schedule of a firm, over a period of time and a given electrical network. With growing importance of environmental impact, some objectives of interest include $\textrm{CO}_2$ emission levels and renewable energy penetration, in addition to the standard generation costs. Some typical constraints include limits on generation levels and up/down times on generation units. This further entails solving a multiobjective mixed integer optimization problem. The related literature has predominantly focused on heuristics (like Genetic Algorithms) for solving larger problem instances. Our major intent in this work is to propose scalable versions of mathematical optimization based approaches that help in speeding up the process of estimating the underlying Pareto frontier. Our contributions are computational and rest on two key embodiments. First, we use the notion of both epsilon constraints and adaptive weights to solve a sequence of single objective optimization problems. Second, to ease the computational burden, we propose a Mccormick-type relaxation for quadratic type constraints that arise due to the resulting formulation types. We test the proposed computational framework on real network data from~\cite{antonio06, multiobjco220} 
and compare the same with standard solvers like Gurobi. Results show a significant reduction in complexity (computational time) when deploying the proposed framework.}
\section{Introduction}
\label{sec:intro}
The standard unit commitment problem~\cite{sc98claus,feltenbasic96}, also known as the OPF (Optimal Power Flow problem) considers the power planning problem of an agent with multiple generators housed in an electrical network. 
Besides cost minimization, this further requires meeting constraints on consumer demand, ramping, switching decisions of generators, the respective up-and-down times~\cite{shahid95cons,ramp08}, and network transmission levels. Up and downtime constraints are common in settings with conventional generators that use coal like fuels, where units cannot be instantly switched on or off due to factors on thermal stress~\cite{updownoren10,birge04}. The switching and up/down time constraints involve discrete variables and this translates the setting into a mixed binary optimization problem in large dimensions~\cite{bertrob13,oren11ucomm}. For incorporating transmission constraints, both DC (Direct Current) type approximations and AC (Alternating Current) type formulations have been discussed in literature~\cite{acopf99,4form12,4form13}.

The last two decades have seen some focus on increasing the mix of renewable energy into these markets. Renewable energy, though environment friendly comes at an additional cost of uncertainty. Besides generation costs, the other important metrics that have been considered so far in literature are carbon emission levels and penetration levels of renewable energy~\cite{multiobjco220,multiobj-ucomm1}. These lead to multiobjective versions of the unit commitment problem. For instance~\cite{multiobj-ucomm1} and ~\cite{multiobjco220} respectively focus on biobjective and triobjective problems with objectives like production costs, maintenance costs, $\textrm{CO}_2$ emissions, and sulphur emissions. These optimization problems do not generally have a single solution, but a family of solutions, usually denoted by a Pareto-frontier. These solutions represent trade-offs between the multiple objectives of interest. Let $x \in \mathbb{R}^{n}$ refer to the space of input variables and $f_j: \mathbb{R}^{n} \rightarrow \mathbb{R}$ represent objective functions for $j = 1,\hdots, m$. Then, for any two such Pareto optimal solutions $x_1^*$ and $x_2^*$ (also called non-dominated points), there exist $i, j \in [m] = \{1,\ldots,m\}$ such that $f_i(x_1^*) < f_i(x_2^*)$ and $f_j(x_1^*) >  f_j(x_2^*)$. 
The quality of Pareto fronts are usually determined by using the hypervolume metric~\cite{audet18,stk12}. More details are listed in the subsequent sections of this work.

In the context of SOUC (Single Objective Unit Commitment), the standard has been to use branch and bound schemes~\cite{mit-bnb,landbnb60}.
Solvers like \texttt{cplex} and \texttt{gurobi} use a host of cutting plane type methods~\cite{mit-bnc2,gomory60,dey07strong,fenchel93boyd} and other practical enhancements (includes heuristics) on top of traditional BnB and several other practical enhancements. Other techniques of interest have been genetic algorithms and heuristics that possibly get good lower bounds on the problem~\cite{ga08}. Some approaches consider Lagrangian based schemes, where the coupling constraints are relaxed~\cite{bertrob13,oren11ucomm}. A detailed test case study of a few heuristics and relaxation based decomposition based for SOUC problems is presented in \cite{watson08}.
Algorithms in the context of MOUC (multiobjective unit commitment) have focused on heuristics and genetic algorithms~\cite{ucomm-ga1,ucomm-ga2,ucomm-ga3,ucomm-ga4,ucomm-ga5,multiobj-ucomm1}. There has been some effort from the standpoint of mathematical programming based approaches~\cite{multiobjco220} for biobjective problems. However, this is a stylized case and centers on a very low dimensional example.

Pareto based avenues to solve multiobjective optimization problems either rely on weights~\cite{Weck06} or epsilon constraints~\cite{epsilon-cons09}. We also specifically note that it has been well accepted in literature that the use of adaptive weights for scalarization~\cite{Weck06} is computationally superior in comparison to uniform weights. 
When the objective functions are nonlinear, using either adaptive weights or epsilon constraints lead to solving mixed-integer optimization problems with nonlinear constraints. 
Our work proposes the use of McCormick type relaxations to alleviate the difficulty in handling the nonlinear constraints. Our contributions are two fold and are stated as follows.
\begin{itemize}
\item \textbf{McCormick relaxations:} Along the lines of~\cite{antonio06}, we consider bi-objective problems (with quadratic objectives) aiming to minimize production costs and carbon emissions. We further note that~\cite{antonio06} considers a single objective version of the problem with ``unit weights'', while we attempt to compute the Pareto frontier. With the introduction of adaptive weights or epsilon constraints, the resulting sequence of problems present quadratic objectives and quadratic constraints in addition to binary variables and standard linear constraints on startup, capacity, and running conditions. While solvers like \texttt{gurobi} handle quadratic constraints, they come at a cost of a huge computational overhead. As the major contribution of this work, we relax such quadratic constraints by means of McCormick relaxations~\cite{mccormick11,tight-mccormick14}. Noting that McCormick relaxations also introduce additional binary variables, we analyze and compare the tradeoffs with solving the quadratically constrained problem. Additionally, it can be observed that the accuracy of our approach depends on the coarseness of McCormick relaxations. 
We compare solutions and computational overhead across different levels of fineness of McCormick relaxations, with the primitive approach having quadratic constraints. We observe that coarser McCormick relaxations lead to good quality Pareto frontiers for comparatively lower computational budget levels.
\item \textbf{Practical use case:} We use the practical datasets from~\cite{antonio06,multiobjco220}. This consists of both thermal and hydro generators. We choose a case with 20 thermal and 10 hydro generators. While we consider only two objectives for this work, we note that the framework extends without loss of generality to larger dimensions (in the space of objectives). We plan to investigate along these lines as part of future research.
\end{itemize}
The rest of the paper is organized into four sections. Section~\ref{sec:model}
discusses about the basics of multiobjective optimization (MOO) and the formulation specifics related to MOUC. Section~\ref{sec:algorithm} proposes the algorithmic framework and the related McCormick relaxations. Section~\ref{sec:numerics} presents numerical results and comparisons on practical test cases from~\cite{matpower11}. We conclude in section~\ref{sec:conclusion}.

\section{Basics and Formulation}
\label{sec:model}
In this section, we introduce the basics on multiobjective optimization (MOO) and present its mathematical formulation. This is ordered into three subsections, where the first two expand on Pareto solutions and MOO algorithms respectively, while the third portion deals with details of the problem formulation.
\subsection{Pareto Optimality}

MOO deals with several conflicting objectives simultaneously, leading to a set of solutions with different trade-offs, known as Pareto optimal solutions or non-dominated solutions. In MOO, the goal is to find a balance between conflicting objectives rather than a single optimal solution. This involves two tasks. The former refers to finding Pareto optimal solutions through an optimization process, while the latter points to selecting the most preferred solution based on preferences of the decision-maker (DM). A DM can express preference information, assuming that less is preferred to more for each objective \cite{miettinen1999nonlinear}. MOO problems can be mathematically represented as follows~\cite{branke2008multiobjective,Weck06,pike2019multi}.
%

\[
\begin{aligned}
& \text{minimize} \quad F(x) = \left[f_1(x), \hspace{0.1mm}f_2(x),\hspace{0.1mm}  \ldots, \hspace{0.1mm} f_k(x)\right] \\
& \text{subject to} \quad x \in S,  \\
\end{aligned}
\]

where \(k \geq 2\) and objective functions \(f_i: \mathbb{R}^n \to \mathbb{R}\) are to be minimized simultaneously. The decision variables are denoted by $x \in S \subseteq \mathbb{R}^{n}$, where $S$ refers to the set of constraints (assumed to be closed and convex throughout the discussion). 
In MOO, a solution is considered optimal if no component of the objective vector can be improved without worsening at least one other components. A decision vector \(x^* \in S\) is \textbf{Pareto optimal} if there is no other \(x\) in the feasible region such that \(f_i(x) \leq f_i(x^*)\) for all \(i = 1, \ldots, k\) and \(f_j(x) < f_j(x^*)\) for at least one index \(j\). We denote the set of Pareto optimal decision vectors by \(P(S)\), and the corresponding set of Pareto optimal objective vectors is \(P(F)\). We formally define the concept of dominated points and non-dominated solutions as follows.
\begin{definition}
A point $x^a$ is said to dominate point $x^b$ iff $f_j(x^{a}) < f_j(x^b)$ for at least one index $j$ and $f_i(x^a) \leq f_i(x^b)$, $\forall \hspace{1mm} i \neq j$. That is,
$$ x^a \succcurlyeq x^b \hspace{2mm} \textit{iff} \hspace{2mm} f_j(x^a) < f_j(x^b) \hspace{2mm} \text{and} \hspace{2mm} f_i(x^a) \leq f_i(x^b).$$
A point $x^c$ is non-dominated iff there exists $j \in \left\{1,\hdots,k \right\}$ such that $f_j(x^{c}) < f_j(x)$ for at least one index $j$ and $\forall x \in S$.
\end{definition}
The non-dominated set of points constitute the Pareto frontier.
To estimate the quality of Pareto frontiers, besides indicators deployed by literature (covered later), the \texttt{utopia} and \texttt{nadir} points are particularly useful.
The \textbf{ideal objective vector} \(F_{\textrm{utopia}}\in \mathbb{R}^k\) represents the lower bounds, obtained by minimizing each objective function individually. That is, $F_{\textrm{utopia}} = \left[f_{1}(x^{1,*}), \hdots, f_{k}(x^{k,*}) \right]$, where $x^{i,*} = \textrm{argmin}_{x \in S} f_i(x)$. Similarly, the \texttt{nadir} point refers to the worst case of all the objectives. That is, $F_{\textrm{nadir}} = \left[f_{1}(\hat{x}^{1}), \hdots, f_{k}(\hat{x}^{k}) \right]$, where $\hat{x}^{i} = \textrm{argmax}_{x \in S} f_i(x)$.

\paragraph{Hypervolume Indicator:} The hypervolume indicator is a widely-used performance metric in multiobjective optimization. This metric measures the volume of the objective space dominated by a set of solutions relative to a reference point \cite{audet18,stk12,beume2006faster}. The hypervolume metric provides a single scalar value that reflects both how close the solutions are to the optimal set and how well they are spread out across the objective space simultaneously. Mathematically, given a set of non-dominated solutions \( S \) and a reference point \( r \), the hypervolume \( HV \) is defined as
\[
HV(S, r) = \lambda_{n} \left( \bigcup_{s \in S} [s; r] \right)
\]
where \( \lambda_{n} \) denotes the $n$-dimensional Lebesgue measure. The reference point is usually chosen to be dominated by all solutions in the Pareto front. A common practice is to set the reference point slightly worse than the worst value observed in the objective functions across all solutions. While there are other metrics that have been used in literature for quantifying the quality of Pareto solutions, we note that this is beyond the scope of our work. The reader is asked to refer to~\cite{audet18} for more details.
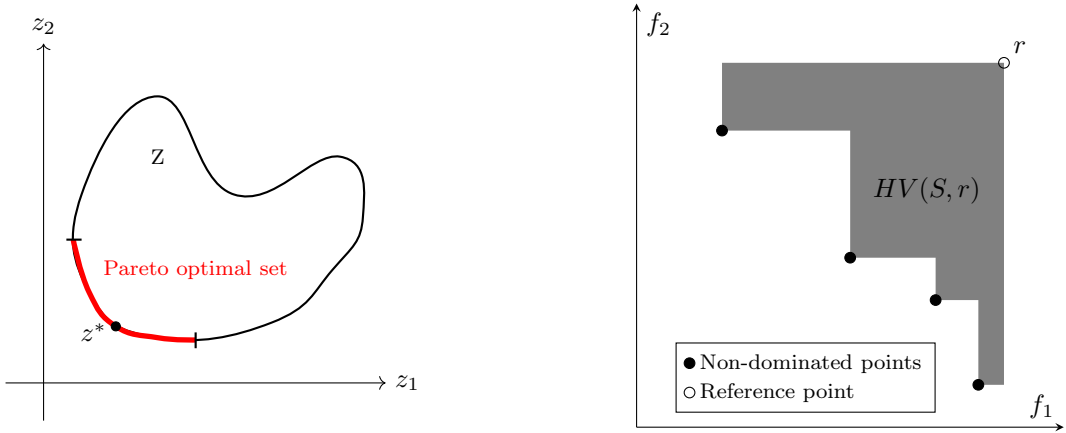
\begin{figure}[h]
\begin{minipage}{0.5\textwidth}
 \centering
  \begin{tikzpicture}
   \draw[->] (-0.5,0) -- (4.5,0) node[right] {$z_1$};
   \draw[->] (0,-0.5) -- (0,4.5) node[above] {$z_2$};
   \draw[thick] plot[smooth cycle, tension=0.8] coordinates{(0.6,1.2) (0.5,2.5) (1.5,3.8) (2.5,2.5) (3.9,3) (4.2,2.2) (3.8,1.5) (3,0.8) (1.5,0.6)};
   \draw[line width=0.7mm, red] plot[smooth, tension=0.8] coordinates{(0.39,1.9) (0.6,1.2) (0.95,0.75) (1.5,0.6) (2,0.57)};
   \draw[thick] (0.3,1.9) -- (0.5,1.9);
   \draw[thick] (2,0.46) -- (2,0.67);
   \node[red, font=\footnotesize] at (2,1.5) {Pareto optimal set};
   \node[black, font=\footnotesize] at (1.5,3) {Z};
   \filldraw[black] (0.95,0.75) circle (1.7pt) node[left, yshift=-0.5ex] {$z^*$};
 \end{tikzpicture}
\end{minipage}
\begin{minipage}{0.5\textwidth}
\centering
\begin{tikzpicture}
        \begin{axis}[
            width=7.2cm,
            height=7.2cm,
            xlabel={$f_1$},
            ylabel={$f_2$},
            xmin=0, xmax=5,
            ymin=0, ymax=5,
            axis lines=middle,
            xtick = \empty,
            ytick = \empty,
            grid=none,
            legend style={at={(0.7,0.2)}, anchor=north east, font=\footnotesize},
            legend cell align={left}
        ]
            \fill[gray, opacity=0.2]
            (axis cs:4.3,0.5) --
            (axis cs:4,0.5) --
            (axis cs:4,1.5) --
            (axis cs:3.5,1.5) --
            (axis cs:3.5,2) --
            (axis cs:2.5,2) --
            (axis cs:2.5,3.5) --
            (axis cs:1,3.5) --
            (axis cs:1,4.3) --
            (axis cs:4.3,4.3) --
            cycle;
            \addplot[only marks, mark=*, mark options={scale=1}]
            coordinates {(1,3.5) (2.5,2) (3.5,1.5) (4,0.5)};
            \addplot[only marks, mark=o, mark options={scale=1}]
            coordinates {(4.3,4.3)};
            \legend{Non-dominated points, Reference point}
            \node at (axis cs:3.4,2.8) {$HV(S, r)$};
            \node at (axis cs:4.3,4.3) [above right] {$r$};
        \end{axis}
    \end{tikzpicture}
\end{minipage}
\caption{Feasible objective region and Pareto optimal set (Left) and the hypervolume indicator for a bi-objective problem. Plots are self generated.}
\end{figure}

\subsection{MOO Algorithms}
In this sub-section, we will discuss about two classic Pareto based methods to solve multiobjective optimization problems. We note that there exist other non-Pareto type methods for MOO like lexicographic techniques. The focus of this paper is however restricted to only Pareto methods. 

\subsubsection{Weighted Sum Method}
This technique focuses on scalarizing the set of objectives into one objective~\cite{branke2008multiobjective,deb2001multi,aswin2021}. The process is repeated for multiple choices of such ``weights'', which further results in solving a sequence/series of single objective optimization problems.
\begin{equation}
\begin{aligned}
& \min_{x} \hspace{1mm}\sum_{i=1}^{k} w_i f_i(x) \\
& \text{subject to} \hspace{1mm} x \in S,  \\
\label{weight}
\end{aligned}
\end{equation}
where \( w_i \geq 0 \) for all \( i \in \left\{1, \ldots, k \right\}\) denotes the set of weights and \( \sum_{i=1}^{k} w_i = 1 \). Weights can be chosen both in uniform and adaptive ways. Say, as an example for the biobjective case, $w_1 = j/10$ and $w_2 = 1-w_1$ for $j = \left\{ 0,1,\hdots,10\right\}$ represents a case with eleven sets of uniform weights. A major disadvantage with uniform weights is that the Pareto frontier may not be uniformly spaced in terms of function values. Recovery of a good Pareto frontier can require many sets of weights and hence can turn out to be computationally expensive. Alternative approaches deploy the concept of adaptive weights~\cite{Weck06,aswin2021}, where weights are sequentially refined based on evaluated objective values. Our framework will focus on adaptive weights and we will discuss more details in the subsequent sections. Pareto frontiers corresponding to MOO problems can be both convex or nonconvex, irrespective of the convexity properties of the objective functions. In nonconvex settings, weighted schemes are not necessarily guaranteed to obtain all points on the Pareto frontier~\cite{das1997closer}. In these problem instances, the method with epsilon-constraints has proven to be better (to be discussed next).
\subsubsection{Epsilon Constraints Method}
In this scheme, a sequence of similar single objective optimization problems are solved. However, the objectives are not scalarized using weights. Instead, one objective function is chosen for optimization, while constraints are imposed on the other objectives. This formulation was first presented in~\cite{haimes1971bicriterion} and can be stated as follows.
\[
\begin{aligned}
& \text{minimize}\hspace{2mm} \quad f_l(x) \\
& \text{subject to} \quad  f_j(x) \leq \varepsilon_j \quad \text{for all} \quad j = 1, \ldots, k, \quad j \neq i,  \\
& \hspace{21mm} x \in S.
\end{aligned}
\]
Note that \(i \in \{1, \ldots, k\}\) and \(\varepsilon_j\) are upper bounds for the objectives (with \(j \neq i\)). 
A solution \(x^* \in S\) is Pareto optimal if and only if it solves the $\varepsilon$-constraint problem for every \(l = 1, \ldots, k\), where \(\varepsilon_j = f_j(x^*)\) for \(j = 1, \ldots, k\) and \(j \neq i\). A unique solution of the $\varepsilon$-constraint problem is Pareto optimal for any upper bounds, which is proven in \cite{miettinen1999nonlinear}. Therefore, Pareto optimality either involves solving \(k\) different problems (which increases computational cost) or obtaining a unique solution (which can be difficult to verify). A significant advantage of the $\varepsilon$-constraint method is that it does not require the problem to be convex, unlike the weighted sum method. Therefore, it can be used for both convex and non-convex Pareto frontiers \cite{deb2001multi}.
However, choosing the appropriate upper bounds \(\varepsilon_j\) can be complicated, especially as the number of objective functions increases. If the bounds are not chosen within the minimum and maximum values of each objective function, the feasible region would become empty, and no solution would be found.

\subsection{Problem Formulation}
We formally describe the problem of interest in this subsection with all the necessary mathematical details. At a very high level, the unit commitment problem is an optimization task in power system operations, where the goal is to determine the optimal schedule for turning power generation units on and off while meeting electricity demand and adhering to various operational constraints \cite{sc98claus,feltenbasic96,birge04}. Our focus is on a single firm scenario with a mix of both renewable and non-renewable energy sources. These generators operate over various time intervals, reflecting the need for flexibility and precision in power system operations. Non-renewable units like thermal generators are typically more expensive to operate and emit higher levels of pollutants, but they offer stable and controllable power output. On the other hand, hydro and renewable units, while generally more environment friendly and cost-effective, are subject to variability due to water availability and environmental regulations~\cite{hobbs2006next}.
We focus on three primary objectives: minimizing operational costs, reducing emissions, and maximizing renewable energy penetration. Each of these objectives is modeled with specific cost functions and constraints, reflecting the economic, environmental, and sustainability goals of the power system. The operational costs include startup costs, running costs, and generation costs, each modeled to capture the different aspects of power generation economics. Emissions reduction is modeled through quadratic functions of power generation, representing the nonlinear relationship between power output and emissions~\cite{multiobj-ucomm1}.
\paragraph{Setup and Notation:} 
Let the planning horizon (days or weeks) be denoted by \(T\). Production decisions are determined for discrete time periods \(t\) within \(\mathcal{T} = \{1, \ldots, T\}\), where these intervals can span minutes, hours, or days depending on operational needs and grid demands. The firm operates multiple generators, indexed by \(\mathcal{I} = \{1, \ldots, I\}\), where each generator can be in an operational state of either being on or off. Once the decision to activate or deactivate a generator is made, these are constrained by minimum up time and down time due to factors like thermal stresses and the need for equipment to cool down or warm up  \cite{updownoren10,birge04}. The decision variables determine the operational status of each unit within the system, the amount of power generated, and the startup and shutdown activities over the planning horizon. Specifically, for each generator \(i\) and time \(t\), we define:
\begin{itemize}
    \item \(y_{it}\) - a binary variable representing whether unit \(i\) is started (\(y_{it} = 1\)) or not (\(y_{it} = 0\)) at time \(t\).
    \item \(z_{it}\) - a binary variable representing whether unit \(i\) is on (\(z_{it} = 1\)) or off (\(z_{it} = 0\)) at time \(t\).
    \item \(g_{it}\) - a continuous variable representing the amount of power generated by unit \(i\) at time \(t\).
\end{itemize}
\subsubsection{Objective Functions}
We consider objectives that reflect broader goals of economic efficiency and environmental sustainability. These can be broadly classified as follows.
\paragraph{Generation Costs:}
The operational costs for power generation are multifaceted and can be categorized into three primary types \cite{watson08,frangioni2008solving,multiobj-ucomm1,gollmer2000unit}:
\begin{itemize}
    \item \textbf{Startup Costs:} These costs are incurred each time a generator is switched from the off mode to on. This transition often requires a significant amount of fuel and power, especially for thermal units. The startup costs are modeled as a linear function of the binary decision variable \(y_{it}\), which indicates whether unit \(i\) is started (\(y_{it} = 1\)) at time \(t\).
    \item \textbf{Running Costs:} These are the costs associated with keeping a generator in an operational state, regardless of the actual power generation. This includes costs for idling, where the generator is on but not necessarily producing power. The running costs are modeled as a linear function tied to the binary variable \(z_{it}\). 
    \item \textbf{Generation Costs:} These costs are directly related to the amount of electricity generated by each unit. Specifically, generation costs are modeled as a quadratic function of the power output.
\end{itemize}
The total generation cost for the utility is the aggregate of all these individual costs across all units and time periods, formulated as:
\begin{equation}
    f_c(y,z,g) = \sum_{i \in \mathcal{I}} \sum_{t \in \mathcal{T}} \left( f^{y}_{it}(y_{it}) + f^{z}_{it}(z_{it}) + f^{g}_{it}(g_{it}) \right),
\end{equation}
where \(f^{y}_{it}(.)\), \(f^{z}_{it}(.)\), and \(f^{g}_{it}(.)\) represent the cost functions for startup, running, and generation, respectively.
\paragraph{Emissions:}
Regulatory requirements on environmental aspects require emissions from $\textrm{CO}_2$ and sulfur to be maintained at lower levels. These emissions are particularly more when generating power from coal and oil based generators. In sync with previous literature, we model the emission levels (Tons/Kgs of $\textrm{CO}_2$) by means of quadratic functions. 
These costs are encapsulated using coefficients $\alpha>0$ and $\beta>0$, which vary based on the fuel type and technology of the generator. The mathematical representation of the emissions objective is:
\begin{equation}
f_e = \sum_{i \in \mathcal{I}} \sum_{t \in \mathcal{T}} \left( \beta_{it} g_{it} + \gamma_{it} g_{it}^2 \right).
\end{equation}
\paragraph{Renewable Energy Penetration:}
This objective aims to maximize the contribution of renewable sources like wind and solar relative to conventional fossil-fuel-based generation. This is critical for transitioning towards a more sustainable and less carbon-intensive power grid. The goal can be expressed either as a direct maximization of renewable output or as an optimization of the ratio of renewable to non-renewable generation. The formulations for these objectives are as follows.
\begin{align}
f_r = -\sum_{t \in \mathcal{T}} \left( \sum_{i \in \mathcal{I}_c} g_{it} - \sum_{i \in \mathcal{I}_r} g_{it} \right)
\hspace{2mm} \text{or alternatively,} \hspace{2mm}
f_r = \frac{\sum_{t \in \mathcal{T}} \sum_{i \in \mathcal{I}_r} g_{it}}{\sum_{t \in \mathcal{T}} \sum_{i \in \mathcal{I}_c} g_{it}}.
\end{align}
Note that $\mathcal{I}_r$ and $\mathcal{I}_c$ denote the sets of renewable and conventional generators, respectively.

\subsubsection{Constraints}
The unit commitment problem involves a variety of constraints that ensure the feasibility and reliability of power generation schedules. These constraints cover aspects such as startup and shutdown protocols, generation limits, and the need to meet electricity demand at all times. Below, we detail the key constraints incorporated into our model. As a simple example, lets take the case of coal generators. Once these coal units are switched on, they need to be running for a certain amount of time before being turned off. On a similar note, units cannot be turned on immediately after they are turned off. These are collectively modeled by means of startup and running variables.
\paragraph{Startup and Shutdown Constraints:}
The relationship between startup variables \( y \) and running variables \( z \) is crucial for accurately modeling the operational dynamics of power units \cite{watson08,ramp08}. If a unit is off at time \( t \) (i.e. \(z_{it} = 0\)) and switched on at time \( t+1 \) (i.e. \(z_{it} = 1\)), then \( y_{i,t+1} \) is 1. In all other scenarios, \( y_{i,t+1} \) is 0. The relationship between startup and running states can be expressed as:
\begin{equation}
y_{it} = \max \left( z_{it} - z_{i,t-1}, 0 \right).
\label{Startup1}
\end{equation}
Given our assumptions on \( f^{y}(.) \) being convex, increasing, and positive,
(\ref{Startup1}) can be formulated as follows.
\begin{align}
y_{it} \geq z_{it} - z_{i,t-1}, \hspace{2mm} y_{it} \geq 0.
\end{align}

\paragraph{Up and Down Time Constraints:}
As mentioned earlier, once a unit is turned on, it must remain on for a few time periods (up-time), and similarly, once turned off, it must remain off for some time (down-time). Let the respective ``minimum'' up and down times be denoted by \( L_{i} \) and \( l_{i} \). The related constraints can be modeled as follows. For more details on such models, the reader is asked to refer to~\cite{gollmer2000unit,dentcheva1997solving}. 
\begin{align}
z_{it} - z_{i,t-1} &\leq z_{i,\tau}, \quad \forall t \in \mathcal{T}, \forall \tau \in \{t+1, \ldots, \min(t+L_{i}-1,T)\}, \\
z_{i,t-1} - z_{it} &\leq 1 - z_{i,\tau}, \quad \forall t \in \mathcal{T}, \forall \tau \in \{t+1, \ldots, \min(t+l_{i}-1,T)\}.
\end{align}
\paragraph{Generation Constraints:}
Each generator has specific operational limits defined by its minimum and maximum generation levels, denoted by \( q_{i} \) and \( Q_{i} \), respectively. A unit can only produce power when it is on, thus the generation levels must be bounded by these limits and can be expressed as follows.
\begin{equation}
q_{i} z_{it} \leq g_{it} \leq Q_{i} z_{it}, \quad \forall i \in \mathcal{I}, \forall t \in \mathcal{T}.
\end{equation}
In this work, we do not consider ramp constraints. Note that these constraints can be easily embedded into our model without loss of generality.
\paragraph{Demand Constraints:}
In addition to individual generator limits, the total power generated at each time step must meet or exceed the electricity demand. This requirement is crucial for maintaining grid stability and is represented by the following constraint.
\begin{equation}
\sum_{i=1}^{I} g_{it} \geq d_{t}, \quad \forall t \in \mathcal{T}.
\end{equation}

\subsubsection{Mathematical representation in matrix form}
For the ease of parsing, we transform our formulation into a matrix notation. The MOUC hence can be put forth as follows. 
\begin{align}
\min_{y, z, g} \quad & F(y, z, g) = \left[ f_c(y, z, g), f_e(g), f_r(g) \right] \nonumber \\
\text{subject to} \quad & \{ y, z, g \} \in \mathcal{X},
\label{eq:initprob}
\end{align}
where \( F(y, z, g) \) encompasses the cost, emission, and renewable energy penetration objectives. The feasible region $\mathcal{X}$ is defined by the intersection of the startup, shutdown, up and down time constraints, generation capacities, and demand constraints.
\begin{align}
\mathcal{X} = \left\{ x_b, g \mid A_b x_b + B g \leq 0, \; D g \leq u, \; E x_b \leq v, \; x_b \in \{0, 1\} \right\},
\end{align}
where the matrices \( A_b \), \( B \), \( D \), and \( E \) encode the relationships and bounds of the variables. The specific elements of these matrices are defined as follows:
\begin{align*}
& x_b = \begin{pmatrix} y \\ z \end{pmatrix}, g = \begin{pmatrix} g_{1} \\ \vdots \\ g_{T} \end{pmatrix},
z = \begin{pmatrix} z_{11} \\ \vdots \\ z_{1T} \\ \vdots \\ z_{IT} \end{pmatrix},
A_b = \begin{pmatrix}
0 & q_{1}I \\ \vdots & \vdots \\ 0 & q_{I}I \\
0 & -Q_{1}I \\ \vdots & \vdots \\ 0 & -Q_{I}I
\end{pmatrix},
B = \begin{pmatrix}
-I & \cdots & 0 \\ \vdots & \ddots & \vdots \\
0 & \cdots & -I \\ I & \cdots & 0 \\ \vdots & \ddots & \vdots \\
0 & \cdots & I
\end{pmatrix}, \\
D =& \begin{pmatrix} I & \cdots & I \end{pmatrix},
u = \begin{pmatrix} d_{1} \\ \vdots \\ d_{T} \end{pmatrix},
E = \begin{pmatrix}
E_{1} \\ E_{2} \\ -E_{2}
\end{pmatrix},
E_{1} = \begin{pmatrix}
-1 & 1 & \cdots & 0 & 0 & -1 & \cdots & 0 \\
\vdots & \vdots & \ddots & \vdots & \vdots & \vdots & \ddots & \vdots \\
0 & \cdots & -1 & 1 & 0 & 0 & \cdots & -1
\end{pmatrix}.
\end{align*}
For the sake of simplicity, the set $\mathcal{X}$ can be further denoted as follows. Note that we use the following notation for the rest of our paper.
\begin{align}
\mathcal{X} = \left\{x| Ax \leq b, x_b \in \left\{ 0,1 \right\} \right\}, \hspace{1mm} x = \pmat{g \\ x_b}.
\end{align} 
\section{Algorithms}
\label{sec:algorithm}
As stated earlier in section~\ref{sec:intro}, the MOUC has received less interest from the standpoint of Mathematical programming based approaches.
One of the key reasons for the same is that the computational complexity of these schemes increases significantly as the problem dimension grows, especially when dealing with both continuous and integer variables. One major challenge is that using uniform weights to generate a reasonable Pareto frontier requires evaluating an enormous number of weight combinations. For a detailed discussion of this complexity, the reader is referred to previous works~\cite{Weck06,kim2005adaptive,aswin2021}. The use of adaptive weights~\cite{kim2005adaptive} on the other hand requires lesser number of weight combinations. However, using adaptive weights leads to some complexity due to the introduction of nonlinear constraints into the problem. To alleviate the issue of nonlinear constraints, we propose using McCormick-type relaxations, which approximate the feasible space by applying linear envelopes. This approach is discussed in more detail in the following subsections.
\subsection{Adaptive Weighted Sum (AWS) Method}
\begin{figure}[h]
    \centering
    \begin{subfigure}[b]{0.45\textwidth}
        \centering
        \begin{tikzpicture}[scale=0.8]
            \draw[thick, domain=0.2:4.1, dashed, variable=\x, blue] plot ({\x}, {1.5+0.25*(\x-2)*(\x-6)}) node[midway, left, xshift=14ex, yshift=8ex] {\parbox{2cm}{unknown \\ Pareto front}};

            \draw[->] (0,0) -- (4.5,0) node[right] {$f_1$};
            \draw[->] (0,0) -- (0,4.5) node[above] {$f_2$};


            \draw[thick] (0.85,3) -- (3.1,0.75) node[midway, right, xshift=0ex, yshift=4ex] {\scriptsize\parbox{3cm}{piecewise linearized \\ Pareto front}};

            \filldraw[black] (0.85,3) circle (2pt) node[left] {$P_2$};
            \filldraw[black] (3.1,0.75) circle (2pt) node[below right] {$P_1$};
        \end{tikzpicture}
        \caption{}
    \end{subfigure}
    \hfill
    \begin{subfigure}[b]{0.45\textwidth}
        \centering
        \begin{tikzpicture}[scale=0.8]
            \draw[thick, domain=0.2:4.1, dashed, variable=\x, blue] plot ({\x}, {1.5+0.25*(\x-2)*(\x-6)});

            \draw[->] (0,0) -- (4.5,0) node[right] {$f_1$};
            \draw[->] (0,0) -- (0,4.5) node[above] {$f_2$};

            \draw[thick] (0.85,3) -- (3.1,0.75);

            \filldraw[black] (0.85,3) circle (2pt) node[left] {$P_2$};
            \filldraw[black] (3.1,0.75) circle (2pt) node[below right] {$P_1$};

            \draw[thick, red] (3.1,0.75) -- (2.45,1.4) node[midway, below left] {$\delta_j$};
            \draw[thick, red] (0.85,3) -- (1.5,2.35) node[midway, below left] {$\delta_j$};
            \draw[thick] (3.1,0.75) -- (3.1,1.4) ; 
            \draw[thick] (0.85,3) -- (1.5,3) ; 


            \draw[->, green] (3.1,1.4) -- (2.45,1.4) node[above, xshift=2ex] {$\delta_1$};
            \draw[<-, green] (1.5,2.35) -- (1.5,3) node[right, yshift=-1.5ex] {$\delta_2$};
            \draw[thick, black, dashed] (2.45,1) -- (2.45,2.35) ; 
            \draw[thick, black, dashed] (1,2.35) -- (2.45,2.35) ; 
        \end{tikzpicture}
        \caption{}
    \end{subfigure}
    \vfill
    \begin{subfigure}[b]{0.7\textwidth}
        \centering
        \begin{tikzpicture}[scale=0.8]
            \draw[thick, domain=0.2:4.1, smooth, variable=\x, blue] plot ({\x}, {1.5+0.25*(\x-2)*(\x-6)});

            \draw[->] (0,0) -- (4.5,0) node[right] {$f_1$};
            \draw[->] (0,0) -- (0,4.5) node[above] {$f_2$};

            \filldraw[black] (0.85,3) circle (2pt) node[left] {$P_2$};
            \filldraw[black] (3.1,0.7) circle (2pt) node[below right] {$P_1$};

            \draw[thick, black, dashed] (2.45,0.5) -- (2.45,2.35) node[below, xshift=-2.3ex] {\scriptsize\parbox{3cm}{Feasible \\ region}};
            \draw[thick, black, dashed] (0.5,2.35) -- (2.45,2.35) ; 

            \filldraw[orange] (2.45,1.1) circle (2pt) node[below] {\scriptsize\parbox{3cm}{new \\ solutions}};
            \filldraw[orange] (1.27,2.35) circle (2pt) ;
            \filldraw[orange] (1.82,1.68) circle (2pt);
        \end{tikzpicture}
        \caption{}
    \end{subfigure}
    \caption{The offset distances, $\delta_1$ and $\delta_2$}
\label{fig:algoaws}
\end{figure}
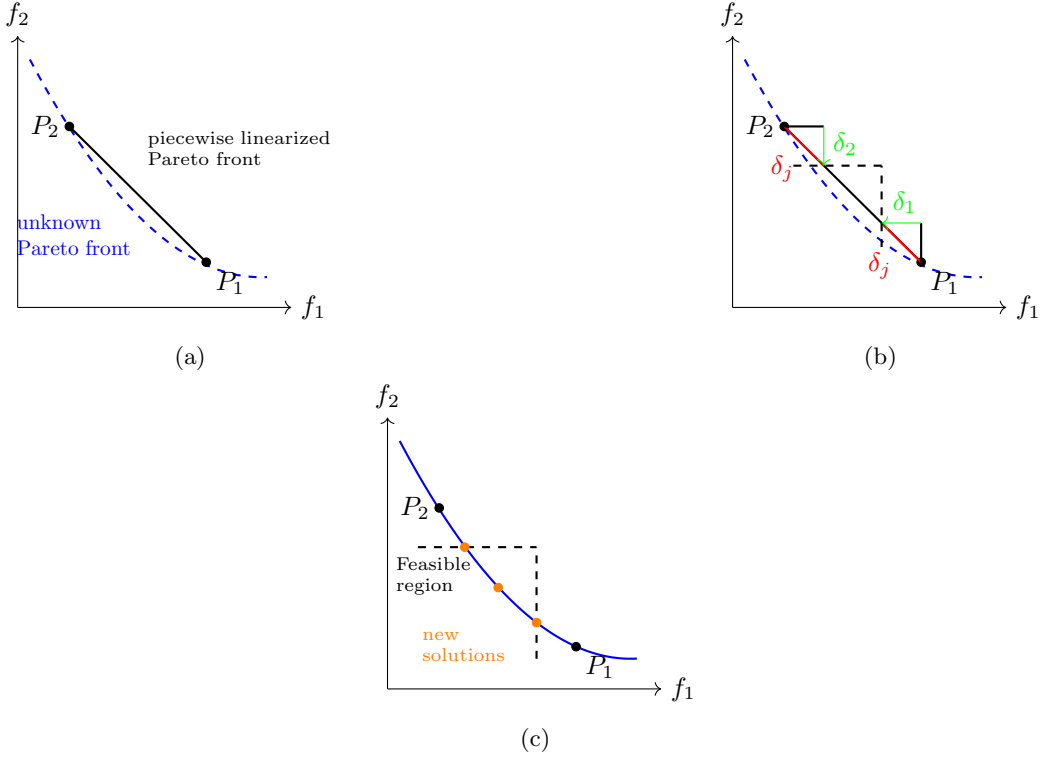

In this work, we focus only on bi-objective problems. The AWS method begins with an initial step similar to the traditional weighted-sum approach but with a large step size for the weighting factor, \( \lambda \). This initial step generates a coarse representation of the Pareto front. The method then identifies regions that need further refinement by calculating distances between neighboring solutions. These identified regions are targeted for sub-optimization with additional inequality constraints. The process is repeated iteratively, refining the Pareto front until a pre-specified resolution is achieved. We summarize the method of adaptive weights from~\cite{kim2005adaptive} for the case of two objectives in figure~\ref{fig:algoaws}. This is subsequently explained stepwise in algorithm~\ref{alg:awsbiobj}.

\begin{algorithm}
\caption{Adaptive Weighted-Sum Method (AWS)}
\label{alg:awsbiobj}
\begin{algorithmic}[1]  
\State \textbf{Input:} Objective functions $f_i$, initial parameters $\delta_J$ (distance threshold), $\varepsilon$ (tolerance), $C$ (scaling factor), $n_{\textrm{initial}}$ (initial number of solutions)
\State \textbf{Output:} Refined Pareto front
\State \textbf{Normalize Objective Functions:}
    \For{each objective function $f_i$}
        \State Compute $f_i^* = \frac{f_i - f_i^{U}}{f_i^{N} - f_i^{U}}$ \Comment{Normalize objective functions}
    \EndFor
    \State \textbf{Initial Weighted-Sum Optimization:}
    \State Perform weighted-sum optimization with step size $\Delta \lambda = \frac{1}{n_{\textrm{initial}}}$.
    \Repeat
        \State \textbf{Segment Length Calculation:}
        \For{each neighboring solution pair}
            \State Calculate Euclidean distance $l_i$ between the two solutions
            \If{$l_i < \varepsilon$}
                \State Remove overlapping solution \Comment{Remove solutions that are too close (overlap)}
            \EndIf
        \EndFor
        \State \textbf{Determination of Further Refinements:}
        \For{each segment $i$}
            \State Compute $n_i = \text{round} \left( C \frac{l_i}{l_{\text{avg}}} \right)$
            \If{$n_i > 1$}
                \State \textbf{Calculation of the Offset Distances:}
                \State Compute $\theta = \tan^{-1} \left( - \frac{P^{y}_1 - P^{y}_2}{P^{x}_1 - P^{x}_2} \right)$
                \State Compute $\delta_1 = \delta_J \cos \theta$, $\delta_2 = \delta_J \sin \theta$

                \State \textbf{Sub-Optimization with Additional Constraints:}
                \State Perform optimization with additional constraints and $n_i$ weighted pairs:
                \[
                \begin{aligned}
                & \min_{x} \hspace{0.5mm} \lambda f_1(x) + (1 - \lambda) f_2(x) \\
                & \text{s.t.} \\
                & f_1(x) \leq P^{x}_1 - \delta_1 \\
                & f_2(x) \leq P^{y}_2 - \delta_2 \\
                & \lambda \in [0, 1]
                \end{aligned}
                \]
            \EndIf
        \EndFor
        \State \textbf{Iteration and Termination:}
        \State Compute lengths of new segments and remove overlapping solutions
    \Until{all segment lengths are below $\delta_J$}
\end{algorithmic}
\end{algorithm}

\paragraph{Stepwise description of the algorithm:} Algorithm~\ref{alg:awsbiobj} can be further elaborated as follows.
\begin{itemize}
\item Initial Weighted-Sum Optimization: This starts with a traditional weighted-sum approach to multiobjective optimization, using a small number of weight combinations \( n_{\text{initial}} \) and a large step size (\( \Delta \lambda \)),
\[
\Delta \lambda = \frac{1}{n_{\textrm{initial}}}.
\]
\item Segment Length Calculation: This step calculates the Euclidean distances between neighboring solutions to identify segments needing refinement. In process, nearly overlapping solutions, where inter-segment lengths less than a prescribed threshold $\epsilon$, are removed.
\item Determination of Further Refinements: The number of further refinements (\( n_i \)) is determined for each segment based on its length relative to the average length of the segment as follows.
\[
n_i = \textrm{round}\left( C \frac{l_i}{l_{\textrm{avg}}} \right),
\]
Note that \( l_i \) is the length of the \( i \)-th segment, \( l_{\textrm{avg}} \) is the average segment length, and \( C \) is a constant. If \( n_i \leq 1 \), there is no further refinement needed for the \( i \)-th segment. Else, the algorithm proceeds with the following step.
\item Calculation of the Offset Distances: For segments needing further refinement, the offset distances (\( \delta_1 \) and \( \delta_2 \)) are calculated from the segment endpoints in the direction of the objective functions. Here, \(\delta_J\) is the prescribed maximum length for segment refinement, and \(P^{x}_1, P^{x}_2, P^{y}_1,\) and \(P^{y}_2\) are the coordinates of the endpoints of the segment in the objective space. The offset distances are calculated using the angle \( \theta \) between the segment and the objective axes as follows.
\[
\theta = \tan^{-1} \left( - \frac{P^{y}_1 - P^{y}_2}{P^{x}_1 - P^{x}_2} \right)
\]
\[
\delta_1 = \delta_J \cos \theta, \quad \delta_2 = \delta_J \sin \theta
\]
\item Sub-Optimization with Additional Constraints: 
Based on the calculated offset distances, additional inequality constraints are imposed and the following optimization subproblem is solved. Note that $\lambda \in [0,1]$.
\[
\begin{aligned}
& \min_{x} \hspace{1mm} & \lambda f_1(x) + (1 - \lambda) f_2(x) \\
& \text{s.t.} & f_1(x) \leq P^{x}_1 - \delta_1 \\
& & f_2(x) \leq P^{y}_2 - \delta_2.
\end{aligned}
\]
Segments that do not converge to optimal solutions are excluded from further refinement, as they correspond to non-convex regions that do not contain Pareto optimal solutions.
\item Iteration and Termination: Calculate the lengths of the new segments and eliminate nearly overlapping solutions. Continue the refinement process until all segment lengths are reduced to below the specified maximum threshold \( \delta_J \).
\end{itemize}
Without loss of generality, this can be extended to more objectives by fixing the weights of all but two objectives. One such generalization is presented in~\cite{Weck06}.

\paragraph{Parameter Selection:}
The Adaptive Weighted-Sum (AWS) method involves several key parameters that must be carefully set to ensure effective optimization. These parameters are the offset distance, the Euclidean distance for determining overlapping solutions, the constant for further refinement (\(C\)), and the number of Pareto front segments in the initial iteration (\(n_{\textrm{initial}}\)). These parameters influence the refinement process and the distribution of solutions along the Pareto front.
The Adaptive Weighted-Sum (AWS) method involves several key parameters that significantly impact optimization performance. These include the offset distance (\(\delta_J\)), the Euclidean distance for identifying overlapping solutions (\(\epsilon\)), the refinement constant ($C$), and the number of initial Pareto front segments \(n_{\textrm{initial}}\). Each of these parameters governs the refinement process and the distribution of solutions along the Pareto front. Succinctly, this can be put forth as follows.
\begin{itemize}
\item \(\delta_J\) must be small to produce denser solutions. It is typically set between 0.05 and 0.2 in the space of normalized objectives.
\item \(\epsilon\) must be smaller than \(\delta_J\) to ensure accurate segmentation.
\item \(C\) being too small prevents further refinement. Large values of \(C\) result in excessive overlapping solutions, further increasing the computational cost.
\item \(n_{\textrm{initial}}\): Here, a balance is needed. Too few divisions may hinder refinement in later stages, while too many can increase computational burden. A typical range of 3 to 10 divisions is recommended.
\end{itemize}


\subsection{McCormick Relaxation}
McCormick envelopes facilitate the computation of solutions to complex optimization problems with bilinear structure by approximating them with linear constraints~\cite{tight-mccormick14,tsoukalas2014multivariate}. 
The use of such relaxations have been common from the standpoint of continuous optimization problems. Moreover, theoretically it has been established in \cite{tsoukalas2014multivariate,bompadre2012convergence} that McCormick relaxation of the product of two functions can achieve quadratic convergence under certain conditions. For our problem of interest, we present Mccormick relaxations in two levels of granularity. 

\subsubsection{Single-Layer McCormick Relaxation}
These form the foundation of the relaxations and these are the weakest version of the approximations. Consider a bilinear term \( w = xy \), where \( x \) and \( y \) are bounded by \( x \in [x_L, x_U] \) and \( y \in [y_L, y_U] \). Along the lines of~\cite{dombrowski2015mccormick}, the relaxation can be constructed using the following linear inequalities.
\begin{enumerate}
    \renewcommand{\labelenumi}{\arabic{enumi}.}

    \item Consider $a = (x - x^L)$ and $\quad b = (y - y^L)$. Since both \( a \) and \( b \) are non-negative:
        \[
        a \cdot b \geq 0 \implies (x - x^L)(y - y^L) = xy - x^L y - xy^L + x^L y^L \geq 0
        \]
        Rearranging, we obtain:
        \[
        w \geq x^L y + x y^L - x^L y^L
        \]

    \item Consider $a = (x^U - x)$ and $b = (y^U - y)$. Again, since \( a \) and \( b \) are non-negative:
        \[
        a \cdot b \geq 0 \implies (x^U - x)(y^U - y) = x^U y^U - x^U y - x y^U + xy \geq 0
        \]
        Rearranging, we obtain:
        \[
        w \geq x^U y + x y^U - x^U y^U
        \]

    \item Considering $a = (x^U - x)$ and $b = (y - y^L)$, we obtain analogously:
        \[
        w \leq x^U y + x y^L - x^U y^L
        \]

    \item Lastly, considering $a = (x - x^L)$ and $b = (y^U - y)$ we obtain:
        \[
        w \leq x y^U + x^L y - x^L y^U
        \]
\end{enumerate}

\begin{figure}[h]
    \centering
    \begin{tikzpicture}
        \draw[->] (-0.5,0) -- (5.5,0) node[right] {$x$};
        \draw[->] (0,-0.5) -- (0,5.5) node[above] {$y$};

        \draw[thick, domain=0.5:4.5, smooth, variable=\x, blue] plot ({\x}, {0.5+0.2*(\x-5)*(\x-5)}) node[above right] {$w = xy$};

        \draw[thick, red, smooth] (0.5,4.55) -- (1.3,0.9) node[pos=0.75, above left, xshift=-2ex] {Underestimators};
        \draw[thick, red, smooth] (1,1.2) -- (4.5,0.55)  ;

        \draw[thick, black, dashed] (0.5,0) -- (0.5,5.5) node[below, yshift=-33.7ex] {$x^L$};
        \draw[thick, black, dashed] (4.5,0) -- (4.5,5.5) node[below, yshift=-33.7ex] {$x^U$};

        \draw[thick, black, dashed] (0,4.55) -- (5.5,4.55) node[left,xshift=-33ex] {$y^U$};
        \draw[thick, black, dashed] (0,0.55) -- (5.5,0.55) node[left,xshift=-33ex] {$y^L$};

        \draw[thick, green, smooth] (0.5,4.55) -- (3.5,3) node[pos=0.75, above right, yshift=1ex] {Overestimators};
        \draw[thick, green, smooth] (3,3.5) -- (4.5,0.55) ;

    \end{tikzpicture}
    \caption{McCormick Envelopes with Underestimators and Overestimators}
    \label{fig:mccormick_envelopes}
\end{figure}
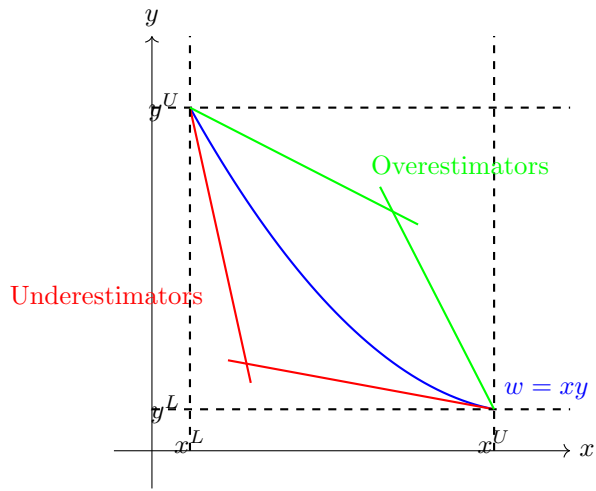

More concisely, \( w = xy \) can further be represented by the \textbf{underestimators}:
\begin{align*}
&w \geq x^L y + x y^L - x^L y^L, \\
&w \geq x^U y + x y^U - x^U y^U
\end{align*}
and the \textbf{overestimators}:
\begin{align*}
&w \leq x^U y + x y^L - x^U y^L, \\
&w \leq x y^U + x^L y - x^L y^U.
\end{align*}
The primary advantage as it can be easily noticed is that the original nonconvex constraints are relaxed by means of convex and linear inequalities. This helps significantly with computational tractability of the problem. The accuracy and tightness of this relaxation remains a question. The next subsection discusses on a more close examination of the same. Our empirical results (discussed later) will cover this in detail.
\subsubsection{Multi-Layer McCormick Relaxation}
These are higher order generalizations, where the domain of the variables is partitioned into multiple segments~\cite{tight-mccormick14}. This yields tighter envelopes. However, such a partitioning also has a disadvantage in that these envelopes are piecewise linear and therefore involve additional binary variables to be introduced into the formulation.  

Consider the bilinear term \( w = xy \) with $x$ and $y$ bounded by \( x \in [x^L, x^U] \) and \( y \in [y^L, y^U] \) as earlier. To apply a multi-layer McCormick relaxation, the domain of \( x \) is partitioned into \( N \) intervals \([x^L_n, x^U_n]\) and the bounds for \([y^L_n, y^U_n]\) are adjusted accordingly. Additionally, binary variables \( q_n \) and a large constant M are introduced into the formulation to accommodate for multiple such piecewise envelopes. Here, the variable \( q_n \) corresponding to the $n$th piece takes the value of ``1'' if the variable ``w'' is in piece (partition) ``\( n \)'' and 0 otherwise. 
The envelopes for partition \( n \) are hence constructed as follows.
\begin{align*}
& w \geq x^L_n y + x y^L_n - x^L_n y^L_n - (1 - q_n) M, \hspace{2mm} \text{Underestimator 1}\\
& w \geq x^U_n y + x y^U_n - x^U_n y^U_n - (1 - q_n) M, \hspace{2mm} \text{Underestimator 2}\\
& w \leq x^U_n y + x y^L_n - x^U_n y^L_n + (1 - q_n) M, \hspace{2mm} \text{Overestimator 1}\\
& w \leq x^L_n y + x y^U_n - x^L_n y^U_n + (1 - q_n) M, \hspace{2mm} \text{Overestimator 2}\\
& q_n \in \{0, 1\}, \hspace{2mm} \sum_{n}^{N} q_n = 1.
\end{align*}
The partitions in the multi-layer McCormick relaxation can be uniform or non-uniform. In the case of uniform partitions, the bounds for $x$ in each partition \( n \) can be defined as follows:
    \begin{align*}
    x_1^L &= x^L, \hspace{2mm} x_N^U = x^U \\
    x_n^L &= x_{n-1}^U, \quad \text{for} \; n = 2, 3, \ldots, N, \\
    x_n^U &= x^L + n \frac{x^U - x^L}{N}, \quad \text{for} \; n = 1, 2, \ldots, N-1.
    \end{align*}
The bounds for \([y^L_n, y^U_n]\) are adjusted accordingly. According to \cite{tight-mccormick14}, multi-layer McCormick relaxation significantly improves the quality of the relaxation and reduces the optimality gap in large-scale optimization problems. This iterative approach allows for better handling of complex constraints and objectives, making it a valuable tool in the context of unit commitment. 

\subsection{Adaptation of Mccormick Relaxations to MOUC}
For completeness, in this section, we expand on the details of the formulation of the Mccormick relaxations, when specifically applied to our algorithm with AWS and epsilon constraints. 
\subsubsection{Adaptive Weights (AWS)} We begin with AWS and present the formulation with Mccormick relaxations for two levels of granularity.
\paragraph{Single Layer Mccormick:}
Consider the matrices \(Q_1\) and \(Q_2\) stated earlier. For the ease of computation, we marginally regularize the matrices by adding a term $\epsilon I$ to both matrices. Note that I refers to an identity matrix and $\epsilon$ was set to be 0.01 times the average value of the elements in \(|Q_1|\) and \(|Q_2|\). Note that $|.|$ denotes the absolute value function. The formulation further involves the following steps.
\begin{itemize}
\item \textbf{Cholesky decomposition}: The regularized matrices \(Q_1\) and \(Q_2\) are factorized using Cholesky decomposition. This results in lower triangular matrices \(L_1\) and \(L_2\), such that $Q_1 = L_1*L_1^{T},\quad Q_2 = L_2*L_2^{T}$.
\item \textbf{Setting variable bounds}: 
The upper bounds \(x_U\) on the decision variables are determined based on the maximum power for continuous variables and are set to ones for binary variables. The lower bounds \(x_L\) are set to zero. For completeness, $x^L \leq x \leq x^U$.
\item \textbf{Auxiliary Variables}: Auxiliary variables \(y_1\), \(y_2\), \(w_1\), and \(w_2\) are introduced such that
\begin{equation}
e^T w_1 = y_1^Ty_1 = \sum_{i=1}^{n} (y_{1,i})^2,
w_{1,i} = y_{1,i}^2, \hspace{1mm} \forall  i=1\dots n,
\end{equation}
\begin{equation}
e^T w_2 = y_2'*y_2 = \sum_{i=1}^{n} (y_{2,i})^2,
w_{2,i} = y_{2,i}^2, \hspace{1mm} \forall i=1\dots n.
\label{eq:w1_w2}
\end{equation}
\end{itemize}
The variables $y_1, y_2$ and the bounds on the same are given as follows.
\begin{align*}
y_1 & = L_1^T x \quad \text{and} \quad y_2 = L_2^T x. \\
y_1^U & = L_1^T x^U, \hspace{1mm} y_1^L = L_1^T x^L, \hspace{1mm}
y_2^U = L_2^T x^U, \hspace{1mm} \text{and} \hspace{1mm} y_2^L = L_2^T x^L.
\end{align*}

Therefore, the constraints marking the refinements can be further expressed as follows. Note that $e$ refers to a vector of ones.
\begin{equation}
e^T w_1 + \textrm{lin}_1^T x \leq P^{x}_1 - \delta_1 \hspace{4mm} \text{and} \hspace{4mm}
e^T w_2 + \textrm{lin}_2^T x \leq P^{y}_2 - \delta_2.
\end{equation}
By relaxing \(w_{1,i} = y_{1,i}^2\) for all  \(i = 1, \dots, n\) and
\(w_{2,i} = y_{2,i}^2 \hspace{1mm} \forall i = 1, \dots, n\) using a single-layer McCormick relaxation, we obtain the following optimization problem with linearized constraints. Note that $\lambda \in \left[ 0,1\right]$ as earlier.
\begin{align*}
 \min_{x} \hspace{1mm} & \lambda f_1(x) + (1 - \lambda) f_2(x) & \\
 \textrm{subject to:} \hspace{1mm} & Ax \leq b, \hspace{1mm} x_b \in \left\{0,1\right\}^{2IT}, \hspace{1mm} x_L \leq x \leq x_U, & (\textrm{original constraints})  \\
& e^T w_1 + \textrm{lin}_1^T x \leq P^{x}_1 - \delta_1, \quad e^T w_2 + \textrm{lin}_2^T x \leq P^{y}_2 - \delta_2, \hspace{3mm} & (\textrm{additional constraints})  \\
& w_1 \geq 2y_1y_1^L - (y_1^L)^2, \hspace{3mm} & (\textrm{underestimator for} \hspace{0.5mm} f_1)\\
& w_1 \geq 2y_1y_1^U - (y_1^U)^2, \hspace{3mm} & (\textrm{underestimator for} \hspace{0.5mm} f_1)\\
& w_1 \leq y_1(y_1^L + y_1^U) - y_1^U y_1^L, \hspace{3mm} &(\textrm{overestimator for} \hspace{0.5mm} f_1)\\
& w_2 \geq 2y_2y_2^L - (y_2^L)^2, \hspace{3mm} & (\textrm{underestimator for} \hspace{0.5mm} f_2)\\
& w_2 \geq 2y_2y_2^U - (y_2^U)^2, \hspace{3mm} &(\textrm{underestimator for} \hspace{0.5mm} f_2)\\
& w_2 \leq y_2(y_2^L + y_2^U) - y_2^U y_2^L, \hspace{3mm} & (\textrm{overestimator for} \hspace{0.5mm} f_2)\\
& y_1 = L_1^T x, \quad y_2 = L_2^T x, \hspace{3mm} & (\textrm{others}).
\end{align*}
Therefore, the linearized constraints can be incorporated into the optimization model, resulting in a new set of constraints \(A_{\textrm{mck1}} \cdot x \leq b_{\textrm{mck1}}\), where:
\begin{footnotesize}
\[
A_{\textrm{mck1}} =
\begin{pmatrix}
A & \mathbf{0} & \mathbf{0} & \mathbf{0} & \mathbf{0} \\
\textrm{lin}^T & \mathbf{0} & \mathbf{0} & \mathbf{1} & \mathbf{0} \\
\textrm{lin2}^T & \mathbf{0} & \mathbf{0} & \mathbf{0} & \mathbf{1} \\
\mathbf{0} & 2\cdot \textrm{diag}(y_1^L) & \mathbf{0} & -\mathbf{I} & \mathbf{0} \\
\mathbf{0} & 2\cdot \textrm{diag}(y_1^U) & \mathbf{0} & -\mathbf{I} & \mathbf{0} \\
\mathbf{0} & \textrm{diag}(-(y_1^U + y_1^L)) & \mathbf{0} & \mathbf{I} & \mathbf{0} \\
\mathbf{0} & \mathbf{0} & 2\cdot \textrm{diag}(y_2^L) & \mathbf{0} & -\mathbf{I} \\
\mathbf{0} & \mathbf{0} & 2\cdot \textrm{diag}(y_2^U) & \mathbf{0} & -\mathbf{I} \\
\mathbf{0} & \mathbf{0} & \textrm{diag}(-(y_2^U + y_2^L)) & \mathbf{0} & \mathbf{I} \\
L_1^T & -\mathbf{I} & \mathbf{0} & \mathbf{0} & \mathbf{0} \\
L_2^T & \mathbf{0} & -\mathbf{I} & \mathbf{0} & \mathbf{0} \\
\end{pmatrix}, \hspace{1mm} b_{\textrm{mck1}} =
\begin{pmatrix}
b \\
P_1^x - \delta_1 \\
P_2^y - \delta_2 \\
y_1^L \cdot y_1^L \\
y_1^U \cdot y_1^U \\
-(y_1^U \cdot y_1^L) \\
y_2^L \cdot y_2^L \\
y_2^U \cdot y_2^U \\
-(y_2^U \cdot y_2^L) \\
\mathbf{0} \\
\mathbf{0} \\
\end{pmatrix}.
\]
\end{footnotesize}
In this optimization problem with linearized constraints, the number of decision variables increases from \( n \) to \( 5n \) due to the inclusion of auxiliary variables \( y_1 \), \( y_2 \), \( w_1 \), and \( w_2 \) alongside the original decision variables \( x \).

\paragraph{Two Layer Mccormick:}
In this case, the bounds for \(y_1\) and \(y_2\) are divided into two pieces as follows.
\[
\begin{aligned}
& y_{1a}^L = y_1^L, \quad y_{1a}^U = \frac{y_1^U + y_1^L}{2}, \quad
y_{1b}^L = y_{1a}^U, \quad y_{1b}^U = y_1^U, \\
& y_{2a}^L = y_2^L, \quad y_{2a}^U = \frac{y_2^U + y_2^L}{2}, \quad y_{2b}^L = y_{2a}^U, \hspace{2mm} \textrm{and} \hspace{3mm} y_{2b}^U = y_2^U.
\end{aligned}
\]
Additionally, we introduce binary variables \(q_a\) and \(q_b\), corresponding to the pieces $a$ and $b$, where \(q_a, q_b \in \{0,1\}\) and \(q_a + q_b = 1\)  \( \forall i = 1, \dots, n\), indicating that only one of the pieces \(a\) or \(b\) is active at a time. Therefore, by relaxing \(w_{1,i} = y_{1,i}^2\) \(\forall i = 1, \dots, n\) and
\(w_{2,i} = y_{2,i}^2 \quad \forall i = 1, \dots, n\) using a two-layer McCormick relaxation, we obtain the following optimization problem with linearized constraints.
\begin{align*}
\min_{x} \quad & \lambda f_1(x) + (1 - \lambda) f_2(x) \\
\text{subject to:} & A x \leq b, \hspace{1mm} x_b \in \left\{ 0,1 \right\}^{2IT}, \hspace{1mm} x_L \leq x \leq x_U, \quad & (\textrm{original constraints}) \\
& e^T w_1 + \textrm{lin}_1^T x \leq P^{x}_1 - \delta_1, \quad e^T w_2 + \textrm{lin}_2^T x \leq P^{y}_2 - \delta_2, \quad & \\ 
& w_1 \geq 2y_1y_{1a}^L -(y_{1a}^L)^2-(1-q_a)M, \\
& w_1 \geq 2y_1y_{1a}^U -(y_{1a}^U)^2-(1-q_a)M, \\
& w_1 \leq y_1 (y_{1a}^L + y_{1a}^U) -y_{1a}^U y_{1a}^L+(1-q_a)M, \\
& w_2 \geq 2y_2y_{2a}^L -(y_{2a}^L)^2-(1-q_a)M, \\
& w_2 \geq 2y_2y_{2a}^U -(y_{2a}^U)^2-(1-q_a)M, \\
& w_2 \leq y_2 (y_{2a}^L + y_{2a}^U) -y_{2a}^U y_{2a}^L+(1-q_a)M, \\
& w_1 \geq 2y_1y_{1b}^L -(y_{1b}^L)^2-(1-q_b)M, \\
& w_1 \geq 2y_1y_{1b}^U -(y_{1b}^U)^2-(1-q_b)M, \\
& w_1 \leq y_1 (y_{1b}^L + y_{1b}^U) -y_{1b}^U y_{1b}^L+(1-q_b)M, \\
& w_2 \geq 2y_2y_{2b}^L -(y_{2b}^L)^2-(1-q_b)M, \\
& w_2 \geq 2y_2y_{2b}^U -(y_{2b}^U)^2-(1-q_b)M, \\
& w_2 \leq y_2 (y_{2b}^L + y_{2b}^U) -y_{2b}^U y_{2b}^L+(1-q_b)M, \\
& q_a \in \left\{0,1 \right\} , q_b \in \left\{0,1 \right\}, \hspace{2mm} q_a + q_b =1.
\end{align*}
Note that \(M\) is a scalar big enough to deactivate a set of constraints when either \(q_{a,i} = 0\) or \(q_{b,i} = 0\). Therefore, this results in a new set of constraints \(A_{\textrm{mck2}} \cdot x \leq b_{\textrm{mck2}}\), where:
\begin{footnotesize}
\[
A_{\textrm{mck2}} =
\begin{pmatrix}
A & \mathbf{0} & \mathbf{0} & \mathbf{0} & \mathbf{0} & \mathbf{0} & \mathbf{0}\\
\textrm{lin}_1^T & \mathbf{0} & \mathbf{0} & \mathbf{1} & \mathbf{0} & \mathbf{0} & \mathbf{0}\\
\textrm{lin}_2^T & \mathbf{0} & \mathbf{0} & \mathbf{0} & \mathbf{1} & \mathbf{0} & \mathbf{0} \\
\mathbf{0} & 2 \cdot \textrm{diag}(y_{1a}^L) & \mathbf{0} & -\mathbf{I} & \mathbf{0} & M \cdot \mathbf{I} & \mathbf{0} \\
\mathbf{0} & 2 \cdot \textrm{diag}(y_{1a}^U) & \mathbf{0} & -\mathbf{I} & \mathbf{0} & M \cdot \mathbf{I} & \mathbf{0} \\
\mathbf{0} & \textrm{diag}(-(y_{1a}^U + y_{1a}^L)) & \mathbf{0} & \mathbf{I} & \mathbf{0} & M \cdot \mathbf{I} & \mathbf{0} \\
\mathbf{0} & \mathbf{0} & 2 \cdot \textrm{diag}(y_{2a}^L) & \mathbf{0} & -\mathbf{I} & M \cdot \mathbf{I} & \mathbf{0} \\
\mathbf{0} & \mathbf{0} & 2 \cdot \textrm{diag}(y_{2a}^U) & \mathbf{0} & -\mathbf{I} & M \cdot \mathbf{I} & \mathbf{0} \\
\mathbf{0} & \mathbf{0} & \textrm{diag}(-(y_{2a}^U + y_{2a}^L)) & \mathbf{0} & \mathbf{I} & M \cdot \mathbf{I} & \mathbf{0} \\
\mathbf{0} & 2 \cdot \textrm{diag}(y_{1b}^L) & \mathbf{0} & -\mathbf{I} & \mathbf{0} & \mathbf{0} & M \cdot \mathbf{I} \\
\mathbf{0} & 2 \cdot \textrm{diag}(y_{1b}^U) & \mathbf{0} & -\mathbf{I} & \mathbf{0} & \mathbf{0} & M \cdot \mathbf{I} \\
\mathbf{0} & \textrm{diag}(-(y_{1b}^U + y_{1b}^L)) & \mathbf{0} & \mathbf{I} & \mathbf{0} & \mathbf{0} & M \cdot \mathbf{I} \\
\mathbf{0} & \mathbf{0} & 2 \cdot \textrm{diag}(y_{2b}^L) & \mathbf{0} & -\mathbf{I} & \mathbf{0} & M \cdot \mathbf{I} \\
\mathbf{0} & \mathbf{0} & 2 \cdot \textrm{diag}(y_{2b}^U) & \mathbf{0} & -\mathbf{I} & \mathbf{0} & M \cdot \mathbf{I} \\
\mathbf{0} & \mathbf{0} & \textrm{diag}(-(y_{2b}^U + y_{2b}^L)) & \mathbf{0} & \mathbf{I} & \mathbf{0} & M \cdot \mathbf{I} \\
L_1^T & -\mathbf{I} & \mathbf{0} & \mathbf{0} & \mathbf{0} & \mathbf{0} & \mathbf{0}\\
L_2^T & \mathbf{0} & -\mathbf{I} & \mathbf{0} & \mathbf{0} & \mathbf{0} & \mathbf{0}\\
\mathbf{0} & \mathbf{0} & \mathbf{0} & \mathbf{0} & \mathbf{0} & \mathbf{I} & \mathbf{I} \\
\end{pmatrix}, \hspace{2mm} b_{\textrm{mck2}} =
\begin{pmatrix}
b \\
P_1^x - \delta_1 \\
P_2^y - \delta_2 \\
(y_{1a}^L)^2 + M \cdot \mathbf{1} \\
(y_{1a}^U)^2 + M \cdot \mathbf{1} \\
-(y_{1a}^U \cdot y_{1a}^L) + M \cdot \mathbf{1} \\
(y_{2a}^L)^2 + M \cdot \mathbf{1} \\
(y_{2a}^U)^2 + M \cdot \mathbf{1} \\
-(y_{2a}^U \cdot y_{2a}^L) + M \cdot \mathbf{1} \\
(y_{1b}^L)^2 + M \cdot \mathbf{1} \\
(y_{1b}^U)^2 + M \cdot \mathbf{1} \\
-(y_{1b}^U \cdot y_{1b}^L) + M \cdot \mathbf{1} \\
(y_{2b}^L)^2 + M \cdot \mathbf{1} \\
(y_{2b}^U)^2 + M \cdot \mathbf{1} \\
-(y_{2b}^U \cdot y_{2b}^L) + M \cdot \mathbf{1} \\
\mathbf{0} \\
\mathbf{1} \\
\end{pmatrix}.
\]
\end{footnotesize}
Here, the number of decision variables is increased from \( n \) to \( 7n \) due to the inclusion of auxiliary variables \( y_1 \), \( y_2 \), \( w_1 \), \( w_2 \), \( q_a \) and \( q_b \) alongside the original decision variables \( x \).

\subsubsection{Epsilon Constraints}
In this section, we just present the case for a single layer Mccormick relaxation. It can be easily observed that this can be generalized to more layers just like the case of adaptive weights without loss of generality. The related $\varepsilon$-constraints in $f_2$ and $f_1$ can be respectively expressed as follows.
\begin{equation}
e^T w_2 + \textrm{lin}_2^T x \leq l_2 + \varepsilon \cdot (u_2 - l_2), \quad e^T w_1 + \textrm{lin}_1^T x \leq l_1 + \varepsilon \cdot (u_1 - l_1)
\end{equation}
By relaxing \( w_{2,i} = y_{2,i}^2 \quad \forall i = 1, \dots, n \) using a single-layer McCormick relaxation, the optimization problem in constraining $f_2$ can be expressed as follows.
\begin{align*}
\min_{x}  \quad & f_1(x) & \\
\textrm{subject to:} \quad & A x \leq b, \hspace{1mm} x_b \in \left\{0,1\right\}^{2IT}, \hspace{1mm} x_L \leq x \leq x_U, \quad & (\textrm{original constraints})\\
& e^T w_2 + \textrm{lin}_2^T x \leq l_2 + \varepsilon \cdot (u_2 - l_2), &\\
& w_2 \geq 2 y_2 y_2^L - (y_2^L)^2, &  \\
& w_2 \geq 2 y_2 y_2^U - (y_2^U)^2, & \\
& w_2 \leq y_2 (y_2^L + y_2^U) - y_2^U y_2^L, & \\
& y_2 = L_2^T x.
\end{align*}
Similarly, the new set of constraints \(A_{\varepsilon} \cdot x \leq b_{\varepsilon}\) can be defined as follows.
\[
A_{\varepsilon} =
\begin{pmatrix}
A & \mathbf{0} & \mathbf{0} \\
\textrm{lin}_2^T & \mathbf{0} & \mathbf{1}\\
\mathbf{0} & 2\cdot \textrm{diag}(y_2^L) & -\mathbf{I} \\
\mathbf{0} & 2\cdot \textrm{diag}(y_2^U) & -\mathbf{I} \\
\mathbf{0} & \textrm{diag}(-(y_2^U + y_2^L)) & \mathbf{I} \\
L_2^T & -\mathbf{I} & \mathbf{0} \\
\end{pmatrix}, \quad b_{\varepsilon} =
\begin{pmatrix}
b \\
l_2 + \varepsilon(u_2-l_2) \\
y_2^L \cdot y_2^L \\
y_2^U \cdot y_2^U \\
-(y_2^U \cdot y_2^L) \\
\mathbf{0} \\
\end{pmatrix}.
\]
Here, the number of decision variables increases from \( n \) to \( 3n \) due to the inclusion of auxiliary variables \( y_2 \) and \( w_2 \) alongside the original decision variables \( x \). The procedure follows analogously for constraining objective $f_1$. 

\section{Numerics}
\label{sec:numerics}
This section presents a comprehensive exploration of optimization techniques applied to the Multi-Objective Unit Commitment (MOUC) problem, focusing on methods that use realistic data for both thermal and hydroelectric units. The techniques evaluated include the Adaptive Weighted Sum (AWS) method and the $\varepsilon$-constraints approach, with a comparison of Gurobi’s quadratic constraint handling and McCormick relaxation techniques. We start by outlining the data sources that were deployed. Next, we discuss the implementation of both uniform and adaptive weights for solving the unit commitment problem.
Lastly, we investigate the impact of varying $\varepsilon$-constraints on the optimization results. By modifying the bounds of these constraints, we assess how these changes influence the optimal solutions, providing insight into the behavior of the objective functions under different constraint conditions.
Hypervolume indicators are used to compare the quality of the Pareto fronts generated by each method. For further specifics on the implementation details, please check the gitHub repository\footnote{Available at: \url{https://scm.cms.hu-berlin.de/aswinkannan1987/unitcommitment-springer-or}}.

\subsection{Data Sources and Overview}
For our study, we utilized data from two widely recognized sources \cite{antonio06, multiobjco220}. The first dataset addresses the ramp-constrained hydro-thermal Unit Commitment problem~\cite{antonio06}, while the second, sourced from~\cite{multiobjco220}, focuses on emissions in the context of low-carbon MOUC. We integrate these datasets into a multiobjective optimization model that simultaneously considers economic cost, CO2 emissions, and sulfur emissions. This approach enables us to evaluate the environmental impact of different power generation schedules and identify strategies that balance emission reduction with economic efficiency. The details of both datasets are provided below for completeness.
\subsubsection{Hydro-thermal dataset}
This dataset includes randomly generated instances that emulate realistic scenarios and has been utilized in several works to validate algorithmic approaches \cite{frangioni2008tighter,frangioni2011sequential,frangioni2016approximated}. The dataset includes detailed information on thermal and hydro units, such as their characteristics, constraints, and operational parameters. Key data includes the planning horizon length (24 periods), the number of thermal units (20), the number of hydro units (10), and a time-series of power demand values (loads) for each interval.
\paragraph{Thermal Units:} Each thermal unit is defined by a unique identifier, cost coefficients (quadratic, linear, and constant), power output limits, minimum up-times and down-times, and start-up cost parameters, including cold and hot start-up costs, a time constant, maximum start-up time, and a fixed start-up cost.

\paragraph{Hydro Unit Description:}
Each hydro unit is defined by a unique identifier, a volume-to-power conversion coefficient, and minimum and maximum flood levels in the basin. The power output of the unit is directly related to the available water volume, which is determined by the flood levels. The maximum power output is calculated by multiplying the conversion coefficient by the maximum flood level, while the minimum power output is determined similarly using the minimum flood level. These calculations ensure that the hydro unit’s power output is proportional to the water available for conversion. Together, these parameters establish the operational power range of each hydro unit. These relationships can be mathematically expressed as follows:
\[
\textrm{maxPower} = \textrm{scalingFactor} \cdot \textrm{volumeToPower} \cdot \textrm{maxFlood}
\]
\[
\textrm{minPower} = \textrm{scalingFactor} \cdot \textrm{volumeToPower} \cdot \textrm{minFlood}
\]
To accurately reflect the power generation limits, the available water resources are used in the hydro units' power output calculations. To prevent excessive reliance on hydro units, a scaling factor less than one is applied, reducing their maximum power output in the optimization model. This adjustment accounts for the economic advantage of hydro energy, given the assumption that it incurs no production costs, and helps balance other operational constraints and objectives.
\subsubsection{Emissions Data}
The emissions data in the model includes CO2 and sulfur emission coefficients, which are used to calculate emissions based on each unit's power output. These coefficients are derived from established references and represent typical emission rates for different types of power plants. The sulfur emission coefficients are particularly useful for understanding the trade-offs between economic and operational decisions within the power sector.

In our model, we assumed that hydro units produce zero emissions, reflecting their characterization as a clean energy source. This assumption emphasizes the environmental benefits of hydroelectric power in our unit commitment analysis. We adapted emission coefficients originally designed for 10 thermal units, replicating them for the 20 thermal units in our study. For the sake of clarity, we note that the ideal penetration of hydro power would be 100 percent if they are not constrained by capacity. This limitation forces the presence of thermal units, which further imposes the regular constraints discussed above in addition to emissions cost.

\textbf{Important: In our study throughout (including tabulations and figures), we note that the generation cost is measured in dollars and the emission levels are quantified by Metric Tonnes.} 

\paragraph{Optimization with Gurobi:}All results presented in this study were computed using the Gurobi optimization solver, with a Matlab-based environment (R2023b) for model formulation. All the computation discussed in this section was performed on a MacAir Machine with a dualcore processor, with a speed of 1.6 GhZ and 8 GB of memory (RAM). The Matlab interface supports quadratic objectives, linear and quadratic inequality constraints, linear equalities, bound constraints, and special ordered set constraints (SOCs). Gurobi facilitates the management of multiple competing objectives, offering two key approaches for handling trade-offs as follows.
\begin{itemize}
    \item \textbf{Blended Approach:} Optimizes a weighted combination of individual objectives.
    \item \textbf{Hierarchical (Lexicographic) Approach:} Prioritizes objectives and optimizes in order of priority, ensuring higher-priority objectives are not degraded.
\end{itemize}
In the context of optimization, quadratic constraints are essential for accurately representing problems that involve quadratic terms in their constraints. Gurobi provides a feature called \texttt{QuadCon} to handle these quadratic constraints effectively.

\subsection{Uniform Weights}
This section explains the use of the uniform weights method to solve a multiobjective unit commitment problem. In our case, since we consider just two objectives, the process of constructing the weight set is fairly straightforward. This takes the following form.
$$w^i_1 = \frac{i}{n}, \hspace{2mm} w^{i}_2 = 1-w^i_1.$$
\begin{table}[h!]
    \centering
    \begin{tabular}{ccc}
        \toprule
        Weight Pair & Objective 1 & Objective 2 \\
        \midrule
        $(0.000, 1.000)$ & $3306171.511$ & $11.921$ \\
        $(0.111, 0.889)$ & $2834788.460$ & $121.099$ \\
        $(0.222, 0.778)$ & $2834788.460$ & $121.099$ \\
        $(0.333, 0.667)$ & $2834788.460$ & $121.099$ \\
        $(0.444, 0.556)$ & $2834788.460$ & $121.099$ \\
        $(0.556, 0.444)$ & $2834788.460$ & $121.099$ \\
        $(0.667, 0.333)$ & $2834788.460$ & $121.099$ \\
        $(0.778, 0.222)$ & $2834788.460$ & $121.099$ \\
        $(0.889, 0.111)$ & $2834788.460$ & $121.099$ \\
        $(1.000, 0.000)$ & $2834788.460$ & $121.099$ \\
        \bottomrule
    \end{tabular}
    \caption{Objective values for $n=2160$ and each uniform weight pair}
\label{tab:objective_values_uniform}
\end{table}
The study examined the impact of varying the number of weight pairs on the objective values in a multiobjective unit commitment problem with 2160 decision variables. The results for 10 different uniformly distributed weight pairs are presented in Table \ref{tab:objective_values_uniform}. The weight pairs (1,0) and (0,1) correspond to single-objective optimizations for objectives $f_1$ and $f_2$ respectively, showing the minimum cost and emissions under constraints. The pair
(0,1) produced notably different results, with objective values of
3306171.511 for cost and
11.921 for emissions, compared to the other weight pairs, which yielded values of 2834788.460 for cost and 121.099 for emissions. This difference underscores the importance of further exploring the Pareto front between these extremes to better understand the trade-offs and potential benefits of intermediate weight pairs.
The performance of the uniform weights method was evaluated in the multiobjective unit commitment problem by testing it with different numbers of decision variables and weight pairs, while documenting the runtimes for each configuration. The number of decision variables was varied by changing the number of time steps. The results are summarized in Tables \ref{uniform_tr} and \ref{uniform_gr}. The total runtime includes the time spent on data processing, cost calculations, constructing constraint and objective matrices, and solving the model with Gurobi for each weight pair. The Gurobi runtime specifically refers to the time spent on solving the optimization problem with Gurobi.
\begin{table}[h!]
\centering
\begin{tabular}{|c|c|c|c|c|c|c|c|}
\hline
\diagbox{\textbf{n}}{\textbf{nr. of pairs}} & \textbf{10} & \textbf{20} & \textbf{30} & \textbf{40} & \textbf{50} & \textbf{60} \\ \hline
\textbf{2160} & 7.55 & 13.58 & 20.54 & 25.60 & 33.33 & 39.29 \\ \hline
\textbf{4320} & 20.84 & 29.02 & 37.39 & 43.24 & 50.88 & 59.00 \\ \hline
\textbf{6480} & 58.41 & 71.49 & 75.72 & 85.57 & 99.22 & 105.73 \\ \hline
\textbf{8640} & 144.47 & 150.60 & 163.94 & 172.44 & 184.53 & 193.76 \\ \hline
\textbf{10800} & 264.68 & 296.26 & 306.82 & 327.89 & 375.07 & 382.80 \\ \hline
\end{tabular}
\caption{Total Runtimes (s) for various numbers of decision variables and uniform weight pairs}
\label{uniform_tr}
\end{table}

\begin{table}[h!]
\centering
\begin{tabular}{|c|c|c|c|c|c|c|c|}
\hline
\diagbox{\textbf{n}}{\textbf{nr. of pairs}} & \textbf{10} & \textbf{20} & \textbf{30} & \textbf{40} & \textbf{50} & \textbf{60} \\ \hline
\textbf{2160} & 5.58 & 11.62 & 17.96 & 23.61 & 31.32 & 36.70 \\ \hline
\textbf{4320} & 7.44 & 15.23 & 23.99 & 32.05 & 39.73 & 48.32 \\ \hline
\textbf{6480} & 9.84 & 18.63 & 28.75 & 38.04 & 48.69 & 57.29 \\ \hline
\textbf{8640} & 12.09 & fz23.32 & 33.31 & 44.08 & 58.06 & 72.45 \\ \hline
\textbf{10800} & 19.37 & 39.08 & 57.60 & 74.90 & 99.98 & 119.53 \\ \hline
\end{tabular}
\caption{Gurobi Runtimes (s) for various numbers of decision variables and uniform weight pairs}
\label{uniform_gr}
\end{table}
The total runtime and Gurobi runtime for various configurations of decision variables and weight pairs show distinct behaviours. As expected, both the total runtime and the Gurobi runtime increase with the number of decision variables and the number of weight pairs.
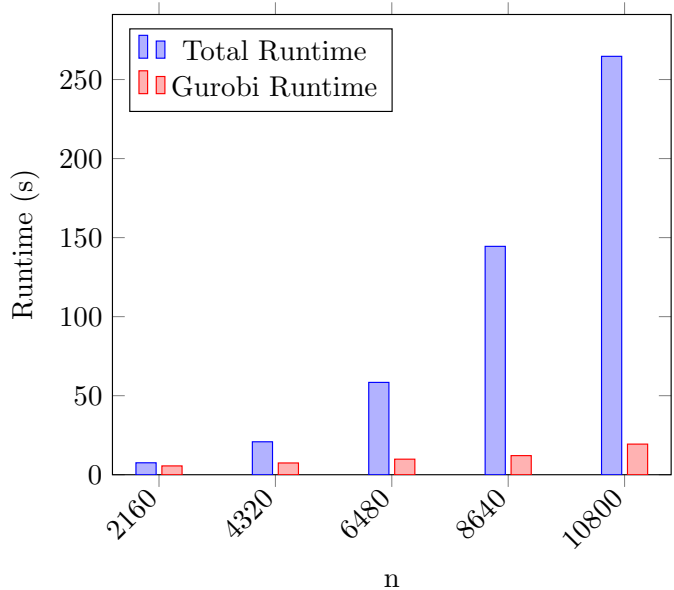
\begin{figure}[h!]
\centering
\resizebox{0.55\textwidth}{!}{ 
\begin{tikzpicture}
    \begin{axis}[
        ybar,
        height=0.3\textheight,
        bar width=7pt,
        xlabel={n},
        ylabel={Runtime (s)},
        symbolic x coords={2160, 4320, 6480, 8640, 10800},
        xtick={2160, 4320, 6480, 8640, 10800},
        x tick label style={rotate=45, anchor=east}, 
        legend pos=north west,
        ymin=0,
        enlarge x limits=0.1
    ]
    \addplot coordinates {(2160,7.55) (4320,20.84) (6480,58.41) (8640,144.47) (10800,264.68)};
    \addplot coordinates {(2160,5.58) (4320,7.44) (6480,9.84) (8640,12.09) (10800,19.37)};
    \legend{Total Runtime, Gurobi Runtime}
    \end{axis}
\end{tikzpicture}
}
\caption{Total and Gurobi Runtimes for 10 weight pairs}
\label{fig:runtime_comparison_k10}
\end{figure}
Figure \ref{fig:runtime_comparison_k10} illustrates the total run times and Gurobi run times for 10 weight pairs. The computational time appears to scale polynomially. This growth in total runtime is primarily due to the presence of increasing overheads as the number of decision variables grows, apart from the optimization time. The Gurobi runtime, which measures the time spent solely on solving the optimization problem, also shows an increase with more decision variables and weight pairs.  However, while the problem size increases, it constitutes smaller proportions of the total runtime compared to other tasks such as data handling and matrix construction.
\subsection{Adaptive Weights}
This section covers the application of our proposed method of adaptive weights, in combination with McCormick relaxations. The computational performance of the adaptive methods is compared against Gurobi's \texttt{QuadCon} and that of weighted methods.
\subsubsection{Formulation}
The adaptive weights method begins with an initial weighted sum optimization to establish a baseline understanding of the objective trade-offs. This initial step uses \( n_{\textrm{initial}} = 10 \) uniformly distributed weight pairs, yielding the results in Table \ref{tab:objective_values_ninitial}.
\begin{table}[h!]
    \centering
    \begin{tabular}{|c|c|c|}
        \hline
        Solutions & Objective 1 & Objective 2 \\ \hline
        $P_1$ & $3306171.511$ & $11.921$ \\ \hline
        $P_2$ & $2834788.460$ & $121.099$ \\ \hline
        $P_3$ & $2834788.460$ & $121.099$ \\ \hline
        $P_4$ & $2834788.460$ & $121.099$ \\ \hline
        $P_5$ & $2834788.460$ & $121.099$ \\ \hline
        $P_6$ & $2834788.460$ & $121.099$ \\ \hline
        $P_7$ & $2834788.460$ & $121.099$ \\ \hline
        $P_8$ & $2834788.460$ & $121.099$ \\ \hline
        $P_9$ & $2834788.460$ & $121.099$ \\ \hline
        $P_{10}$ & $2834788.460$ & $121.099$ \\
        \hline
    \end{tabular}
    \caption{Objective values for $n=2160$ and 10 uniform weight pairs}
    \label{tab:objective_values_ninitial}
\end{table}
After this initial optimization, we proceed through the steps outlined in the previous section.
The results in Table \ref{tab:quadcon_gurobi} show the total and Gurobi runtimes, along with the optimal values for $f_1$ and $f_2$, obtained using the adaptive weights method for different values of decision variables $n$. As $n$ increases (due to more time steps), both runtimes increase significantly, as expected. For example, at $n = 540$, the total runtime is 19.77 seconds, but at $n = 12960$, it rises to 2417.49 seconds. The adaptive weights method generally results in higher runtimes compared to the uniform weights method, due to the additional complexity of dynamically adjusting weights. However, for $n=2160$, the adaptive method successfully explores the Pareto front between two points identified using the uniform weights method. The final Pareto points are stored for further analysis. The majority of the total runtime is spent on the Gurobi optimization process, particularly when using its \texttt{QuadCon} feature to handle quadratic constraints.
\begin{table}
\centering
\begin{tabular}{|c|c|c|c|c|}
\hline
\textbf{n} & \textbf{Total Runtime (s)} & \textbf{Gurobi Runtime (s)} & \textbf{Optimal \(f_1\)} & \textbf{Optimal \(f_2\)} \\
\hline
540 & 19.77 & 19.37 & 717661.658 & 22.246 \\
1080 & 30.02 & 29.44 & 1438587.163 & 41.824 \\
1620 & 63.51 & 62.55 & 2157642.324 & 62.883 \\
2160 & 70.00 & 68.61 & 2878614.292 & 82.412 \\
2700 & 74.96 & 72.83 & 3599634.371 & 101.179 \\
3240 & 77.19 & 73.81 & 4317964.333 & 123.545 \\
3780 & 84.18 & 78.61 & 5037515.019 & 144.279 \\
4320 & 172.57 & 167.66 & 5756169.879 & 166.199 \\
4860 & 214.66 & 207.81 & 6478500.930 & 183.237 \\
5400 & 317.09 & 308.21 & 7199182.782 & 202.468 \\
5940 & 391.10 & 378.19 & 7918177.252 & 223.966 \\
6480 & 466.04 & 479.77 & 8638558.918 & 243.593 \\
7020 & 555.22 & 533.80 & 9359554.651 & 262.397 \\
7560 & 706.52 & 684.14 & 10081211.717 & 280.347 \\
8100 & 718.93 & 695.75 & 10801162.644 & 300.568 \\
8640 & 908.87 & 885.74 & 11520305.713 & 321.839 \\
9180 & 962.22 & 936.23 & 12237578.410 & 345.602 \\
9720 & 1007.20 & 974.50 & 12952620.790 & 372.476 \\
10260 & 1247.08 & 1212.62 & 13677457.939 & 386.039 \\
10800 & 1571.24 & 1531.47 & 14399600.872 & 403.380 \\
11340 & 1853.91 & 1803.21 & 15120516.449 & 422.196 \\
11880 & 2095.26 & 2018.11 & 15837393.454 & 446.522 \\
12420 & 2240.76 & 2176.19 & 16553431.198 & 471.962 \\
12960 & 2417.49 & 2343.15 & 17267926.285 & 499.417 \\
\hline
\end{tabular}
\caption{Total and Gurobi Runtimes for Adaptive Weights using \texttt{QuadCon}}
\label{tab:quadcon_gurobi}
\end{table}
The effectiveness of the adaptive weights method was evaluated using the Hypervolume (HV) metric~\cite{audet18,stk12,beume2006faster}, computed with the \texttt{PyGMO} library\footnote{\url{https://readthedocs.org/projects/pygmo/downloads/pdf/newdocs/}}. The HV was calculated from the non-dominated objective values obtained through the adaptive weights method using \texttt{QuadCon}. For each value of $n$ (number of decision variables), the reference point was set to be approximately 1.2 worse than the worst objective function value across all solutions. The results, presented in Table \ref{tab:HV_aws_quadcon}, show that HV values increase as $n$ increases. However, this increase in HV is not due to better Pareto solutions but because the objective space expands as $n$ grows. 
\begin{table}[h!]
    \centering
    \begin{tabular}{cccccc}
        \toprule
        \textbf{n} & \textbf{Pareto Points} & \textbf{Total (s)} & \textbf{Gurobi (s)} & \textbf{HV ($\cdot 10^7$)} \\
        \midrule
        540 & 6 & 19.77 & 19.37 & 0.22 \\
        1080 & 6 & 30.02 & 29.44 & 0.88 \\
        1620 & 6 & 63.51 & 62.55 & 1.99 \\
        2160 & 6 & 70.00 & 68.61 & 3.51 \\
        2700 & 6 & 74.96 & 72.83 & 5.56 \\
        \bottomrule
    \end{tabular}
    \caption{Hypervolume Indicator for Adaptive Weights for various numbers of decision variables using \texttt{QuadCon}}
    \label{tab:HV_aws_quadcon}
\end{table}

\subsubsection{Quadratic Constraints with single-layer McCormick}
For our next numerical experiment, we relax the quadratic constraints by using a single-layer McCormick relaxation. Table \ref{tab:aws_mccormick} presents the total runtime, Gurobi runtime, and optimized objective values for various numbers of decision variables using the single-layer McCormick relaxation. The number of decision variables presented in table \ref{tab:aws_mccormick} only accounts for the original decision variables \( x \) (and not the auxiliary variables introduced) and is varied by changing the number of time steps considered.
\begin{table}
\centering
\begin{tabular}{|c|c|c|c|c|}
\hline
\textbf{n} & \textbf{Total Runtime (s)} & \textbf{Gurobi Runtime (s)} & \textbf{Optimal \(f_1\)} & \textbf{Optimal \(f_2\)} \\
\hline
540 & 7.59 & 6.87 & 735203.203 & 8.361 \\
1080 & 8.48 & 7.36 & 1472770.926 & 15.530 \\
1620 & 18.69 & 15.77 & 2204594.811 & 25.593 \\
2160 & 19.26 & 15.86 & 2939702.566 & 34.002 \\
2700 & 22.43 & 17.20 & 3674927.273 & 42.352 \\
3240 & 27.54 & 20.41 & 4410151.981 & 50.702 \\
3780 & 32.85 & 23.07 & 5145376.688 & 59.052 \\
4320 & 38.88 & 24.54 & 5880601.472 & 67.401 \\
4860 & 43.48 & 25.09 & 6615826.182 & 75.751 \\
5400 & 51.96 & 28.34 & 7351050.874 & 84.101 \\
5940 & 58.71 & 29.75 & 8086275.514 & 92.451 \\
6480 & 66.52 & 31.29 & 8821500.273 & 100.800 \\
7020 & 77.08 & 31.90 & 9556725.012 & 109.150 \\
7560 & 93.23 & 40.56 & 10291949.721 & 117.500 \\
8100 & 127.17 & 50.91 & 11027174.431 & 125.850 \\
8640 & 135.06 & 61.70 & 11761157.401 & 134.826 \\
9180 & 163.37 & 67.20 & 12497623.824 & 142.549 \\
9720 & 195.49 & 73.70 & 13232028.710 & 151.313 \\
10260 & 217.08 & 86.96 & 13968073.186 & 159.249 \\
10800 & 255.97 & 95.01 & 14701927.560 & 168.290 \\
11340 & 290.33 & 99.40 & 15437081.978 & 176.675 \\
11880 & 327.06 & 114.55 & 16170873.284 & 185.060 \\
12420 & 382.05 & 114.67 & 16908972.098 & 192.648 \\
12960 & 412.56 & 133.63 & 17644196.814 & 200.998 \\
\hline
\end{tabular}
\caption{Total and Gurobi Runtimes for Adaptive Weights using single-layer McCormick}
\label{tab:aws_mccormick}
\end{table}
The results in table \ref{tab:aws_mccormick} demonstrate that using the single-layer McCormick relaxation significantly reduces the runtime in comparison to using the primitive version of \texttt{Quadcon} and \texttt{Gurobi}. Notably, the computational time appears to scale polynomially as the problem size grows. Along with that, e.g. for $n = 2160$, we observe that the adaptive weights
method using single-layer McCormick relaxation effectively explores the Pareto front between the points $P1$ $(3306171.511, 11.921)$ and
$P2$ $(2834788.460, 121.099)$, identified using initial uniform weights method.
Additionally, it can be noted that the \texttt{Gurobi} runtime constitutes a significant portion of the total runtime, though not as much as with using \texttt{QuadCon}. Similar to the previous section, we compute the HV metric using the \texttt{PyGMO} library.
Table \ref{tab:HV_aws_mccormick} provides the HV values for different values of \( n \) alongside the total and Gurobi runtimes, and the number of Pareto points.
\begin{table}[h!]
    \centering
    \begin{tabular}{cccccc}
        \toprule
        \textbf{n} & \textbf{Pareto Points} & \textbf{Total (s)} & \textbf{Gurobi (s)} & \textbf{HV ($\cdot 10^7$)} \\
        \midrule
        540 & 6 & 7.59 & 6.87 & 0.23 \\
        1080 & 6 & 8.48 & 7.36 & 0.91 \\
        1620 & 6 & 18.69 & 15.77 & 2.05 \\
        2160 & 6 & 19.26 & 15.86 & 3.65 \\
        2700 & 6 & 22.43 & 17.20 & 5.75 \\
        \bottomrule
    \end{tabular}
    \caption{Hypervolume Indicator for Adaptive Weights for various numbers of decision variables using single-layer McCormick}
    \label{tab:HV_aws_mccormick}
\end{table}
To ensure a fair comparison, we use the same reference points as in \ref{tab:HV_aws_quadcon}. Comparing the HV values between tables \ref{tab:HV_aws_quadcon} and \ref{tab:HV_aws_mccormick}, we observe that the single-layer McCormick method generally yields higher HV values than the \texttt{QuadCon} method for the same values of \( n \). 

\subsubsection{Quadratic Constraints with 2-layer McCormick}
Next, the quadratic constraints were relaxed using a two-layer McCormick relaxation. Table \ref{tab:quadcon_mccormick_2} presents the total runtime, Gurobi runtime, and final Pareto points derived from this method for various numbers of decision variables using this approach. 
\begin{table}
\centering
\begin{tabular}{|c|c|c|c|c|}
\hline
\textbf{n} & \textbf{Total Runtime (s)} & \textbf{Gurobi Runtime (s)} & \textbf{Optimal \(f_1\)} & \textbf{Optimal \(f_2\)} \\
\hline
540 & 19.14 & 18.06 & 734051.668 & 8.941 \\
1080 & 25.36 & 23.57 & 1468676.718 & 17.593 \\
1620 & 35.38 & 31.16 & 2202605.370 & 26.596 \\
2160 & 44.21 & 36.73 & 2937180.436 & 35.274 \\
2700 & 57.02 & 46.32 & 3671755.841 & 43.951 \\
3240 & 67.92 & 50.45 & 4406330.843 & 52.628 \\
3780 & 71.63 & 51.60 & 5140906.017 & 61.305 \\
4320 & 175.75 & 143.93 & 5875481.108 & 69.983 \\
4860 & 235.29 & 193.27 & 6610056.335 & 78.660 \\
5400 & 255.28 & 201.73 & 7344631.358 & 87.337 \\
5940 & 320.05 & 253.27 & 8079206.517 & 96.015 \\
6480 & 478.27 & 393.82 & 8813796.018 & 104.685 \\
7020 & 557.79 & 453.59 & 9548371.189 & 113.362 \\
7560 & 643.21 & 502.27 & 10282946.309 & 122.039 \\
8100 & 752.10 & 597.37 & 11017507.075 & 130.724 \\
8640 & 796.29 & 622.05 & 11752082.178 & 139.401 \\
9180 & 932.82 & 678.97 & 12486667.585 & 148.073 \\
9720 & 1057.00 & 718.63 & 13221232.476 & 156.755 \\
10260 & 1106.27 & 752.45 & 13955807.595 & 165.433 \\
10800 & 1207.11 & 785.45 & 14690382.752 & 174.110 \\
11340 & 1393.50 & 867.25 & 15424967.464 & 182.782 \\
11880 & 1832.32 & 1179.51 & 16159533.023 & 191.464 \\
12420 & 1946.83 & 1294.42 & 16894108.174 & 200.142 \\
12960 & 2044.38 & 1364.24 & 17628698.231 & 208.811 \\
\hline
\end{tabular}
\caption{Total and Gurobi Runtimes for Adaptive Weights using 2-layer McCormick}
\label{tab:quadcon_mccormick_2}
\end{table}
The 2-layer McCormick method results in higher runtimes than the 1-layer method due to the larger problem size, but its solutions are slightly closer to those achieved with \texttt{QuadCon}. Increasing the number of layers can improve the accuracy further. Despite introducing additional auxiliary variables in the 2-layer method, its runtime remains generally lower than that of \texttt{QuadCon}. Hypervolume (HV) values, presented in Table \ref{tab:HV_mccormick_2}, show a similarly increasing trend with the number of decision variables $n$, comparable to the 1-layer McCormick method. This suggests that while the two-layer method adds complexity, its performance in exploring the objective space is similar to that of the one-layer method.
\begin{table}[h!]
    \centering
    \begin{tabular}{cccccc}
        \toprule
        \textbf{n} & \textbf{Pareto Points} & \textbf{Total (s)} & \textbf{Gurobi (s)} & \textbf{HV ($\cdot 10^7$)} \\
        \midrule
        540 & 6 & 19.14 & 18.06 & 0.23 \\
        1080 & 6 & 25.36 & 23.57 & 0.91 \\
        1620 & 6 & 35.38 & 31.16 & 2.05 \\
        2160 & 6 & 44.21 & 36.73 & 3.65 \\
        2700 & 6 & 57.02 & 46.32 & 5.75 \\
        \bottomrule
    \end{tabular}
    \caption{Hypervolume Indicator for Adaptive Weights for various numbers of decision variables using 2-layer McCormick}
    \label{tab:HV_mccormick_2}
\end{table}
\subsubsection{\texttt{QuadCon} vs. single- and multi-layer McCormick}
\begin{figure}[H]
\centering
\resizebox{0.65\textwidth}{!}{
\begin{tikzpicture}
    \begin{axis}[
        ybar,
        width=\textwidth,
        height=0.5\textheight,
        bar width=4pt,
        xlabel={n},
        ylabel={Gurobi Runtime (s)},
        symbolic x coords={1080, 2160, 3240, 4320, 5400, 6480, 7560, 8640, 9720, 10800, 11880, 12960},
        xtick={1080, 2160, 3240, 4320, 5400, 6480, 7560, 8640, 9720, 10800, 11880, 12960},
        x tick label style={rotate=45, anchor=east},
        legend pos=north west,
        ymin=0,
        enlarge x limits=0.1
    ]
    \addplot coordinates {(1080,29.44) (2160,68.61) (3240,73.81) (4320,167.66) (5400,308.21) (6480,479.77) (7560,684.14) (8640,885.74) (9720,974.50) (10800,1531.47) (11880,2018.11) (12960,2343.15)};
    \addplot coordinates {(1080,7.36) (2160,15.86) (3240,20.41) (4320,24.54) (5400,28.34) (6480,31.29) (7560,40.56) (8640,61.70) (9720,73.70) (10800,95.01) (11880,114.55) (12960,133.63)};
    \addplot coordinates {(1080,23.57) (2160,36.73) (3240,50.45) (4320,143.93) (5400,201.73) (6480,393.82) (7560,502.27) (8640,622.05) (9720,718.63) (10800,785.45) (11880,1179.51) (12960,1364.24)};
    \legend{\texttt{QuadCon}, Single-layer McCormick, 2-layer McCormick}
    \end{axis}
\end{tikzpicture}
}
\caption{Comparison of Gurobi Runtimes for Adaptive Weights using \texttt{QuadCon}, 1-layer and 2-layer McCormick}
\label{fig:comparison_runtimes_aws_gurobi}
\end{figure}
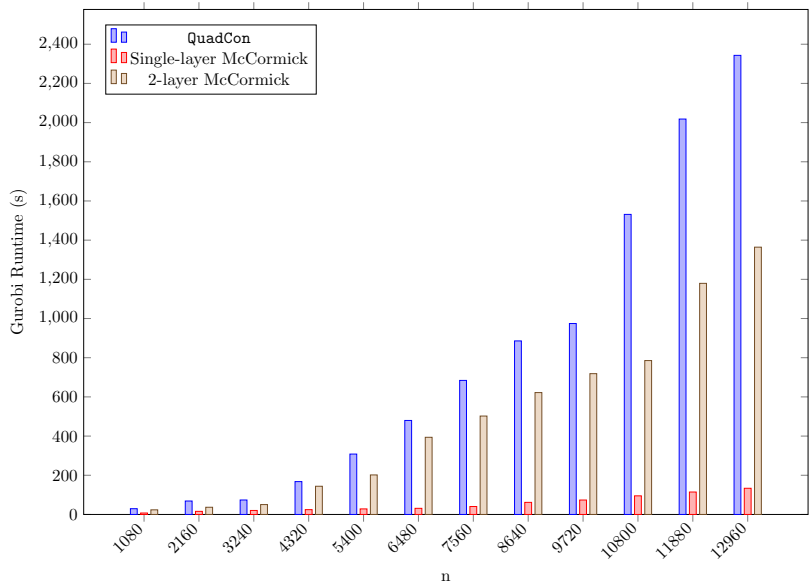

Table \ref{tab:comparison_runtimes_aws_gurobi} and Figure \ref{fig:comparison_runtimes_aws_gurobi} compare the Gurobi runtimes for the adaptive weights method across different values of \( n \) for three distinct approaches to handling quadratic constraints: \texttt{QuadCon}, single-layer McCormick, and 2-layer McCormick relaxations. The results in Table \ref{tab:comparison_runtimes_aws_gurobi} reveal that the single-layer McCormick relaxation consistently achieves the lowest runtime across all values of \( n \), demonstrating its efficiency and applicability to large scale optimization problems.
In contrast, the \texttt{QuadCon} method, which directly handles quadratic constraints, has the highest runtime, especially as the number of decision variables increases. This method's computational intensity is evident from the increase in runtime as \( n \) grows, reaching over 2300 seconds for \( n = 12960 \). The 2-layer McCormick relaxation, while not as fast as the single-layer approach, still offers a considerable improvement over the \texttt{QuadCon} method, thereby finding a balance. For completeness, Figure \ref{fig:comparison_runtimes_aws_gurobi} illustrates the runtime trends for each method as \( n \) increases.
\begin{table}[p]
\centering
\begin{tabular}{|c|c|c|c|}
\hline
\textbf{n} & \textbf{\texttt{QuadCon} (s)} & \textbf{Single-layer McCormick (s)} & \textbf{2-layer McCormick (s)} \\
\hline
540 & 19.37 & 6.87 & 18.06 \\
1080 & 29.44 & 7.36 & 23.57 \\
1620 & 62.55 & 15.77 & 31.16 \\
2160 & 68.61 & 15.86 & 36.73 \\
2700 & 72.83 & 17.20 & 46.32 \\
3240 & 73.81 & 20.41 & 50.45 \\
3780 & 78.61 & 23.07 & 51.60 \\
4320 & 167.66 & 24.54 & 143.93 \\
4860 & 207.81 & 25.09 & 193.27 \\
5400 & 308.21 & 28.34 & 201.73 \\
5940 & 378.19 & 29.75 & 253.27 \\
6480 & 479.77 & 31.29 & 393.82 \\
7020 & 533.80 & 31.90 & 453.59 \\
7560 & 684.14 & 40.56 & 502.27 \\
8100 & 695.75 & 50.91 & 597.37 \\
8640 & 885.74 & 61.70 & 622.05 \\
9180 & 936.23 & 67.20 & 678.97 \\
9720 & 974.50 & 73.70 & 718.63 \\
10260 & 1212.62 & 86.96 & 752.45 \\
10800 & 1531.47 & 95.01 & 785.45 \\
11340 & 1803.21 & 99.40 & 867.25 \\
11880 & 2018.11 & 114.55 & 1179.51 \\
12420 & 2176.19 & 114.67 & 1294.42 \\
12960 & 2343.15 & 133.63 & 1364.24 \\
\hline
\end{tabular}
\caption{Comparison of Gurobi Runtimes for Adaptive Weights using \texttt{QuadCon}, 1-layer and 2-layer McCormick}
\label{tab:comparison_runtimes_aws_gurobi}
\end{table}

\begin{figure}[H]
\centering
\resizebox{0.65\textwidth}{!}{
\begin{tikzpicture}
    \begin{axis}[
        ybar,
        width=\textwidth,
        height=0.5\textheight,
        bar width=4pt,
        xlabel={n},
        ylabel={Total Runtime (s)},
        symbolic x coords={1080, 2160, 3240, 4320, 5400, 6480, 7560, 8640, 9720, 10800, 11880, 12960},
        xtick={1080, 2160, 3240, 4320, 5400, 6480, 7560, 8640, 9720, 10800, 11880, 12960},
        x tick label style={rotate=45, anchor=east},
        legend pos=north west,
        ymin=0,
        enlarge x limits=0.1
    ]
    \addplot coordinates {(1080,30.02) (2160,70.00) (3240,77.19) (4320,172.57) (5400,317.09) (6480,466.04) (7560,706.52) (8640,908.87) (9720,1007.20) (10800,1571.24) (11880,2095.26) (12960,2417.49)};
    \addplot coordinates {(1080,8.48) (2160,19.26) (3240,27.54) (4320,38.88) (5400,51.96) (6480,66.52) (7560,93.23) (8640,135.06) (9720,195.49) (10800,255.97) (11880,327.06) (12960,412.56)};
    \addplot coordinates {(1080,25.36) (2160,44.21) (3240,67.92) (4320,175.75) (5400,255.28) (6480,478.27) (7560,643.21) (8640,796.29) (9720,1057.00) (10800,1207.11) (11880,1832.32) (12960,2044.38)};
    \legend{\texttt{QuadCon}, Single-layer McCormick, 2-layer McCormick}
    \end{axis}
\end{tikzpicture}
}
\caption{Comparison of Total Runtimes for Adaptive Weights using \texttt{QuadCon}, 1-layer and 2-layer McCormick}
\label{fig:total_runtime_comparison_aws}
\end{figure}
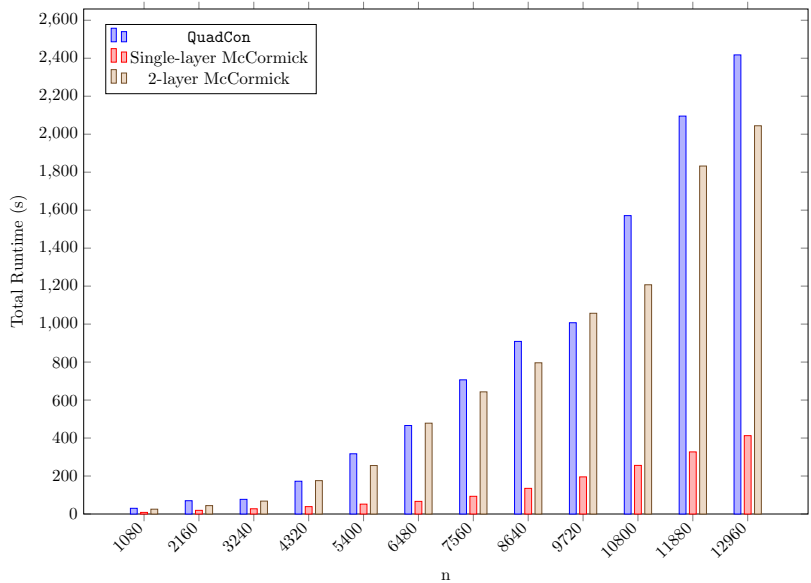

\paragraph{Overhead:} To assess the efficiency of the methods in handling quadratic constraints, both the Gurobi solve time and the additional overhead from auxiliary variables and constraints must be considered. Table \ref{tab:comparison_total_runtimes_aws} and Figure \ref{fig:total_runtime_comparison_aws} present the total runtimes for the \texttt{QuadCon}, single-layer McCormick, and 2-layer McCormick methods across various values of \( n \). The total runtime analysis shows that the single-layer McCormick method is the fastest overall, with lower total runtimes compared to the 2-layer McCormick and \texttt{QuadCon} methods, especially for larger problem sizes. Although the Gurobi solve times for the McCormick methods are faster, their overhead from extra auxiliary variables and larger matrices reduces the performance gains. Future work could focus on optimizing this overhead to improve computational efficiency.
\begin{table}[p]
\centering
\begin{tabular}{|c|c|c|c|}
\hline
\textbf{n} & \textbf{\texttt{QuadCon} (s)} & \textbf{Single-layer McCormick (s)} & \textbf{2-layer McCormick (s)} \\
\hline
540 & 19.77 & 7.59 & 19.14 \\
1080 & 30.02 & 8.48 & 25.36 \\
1620 & 63.51 & 18.69 & 35.38 \\
2160 & 70.00 & 19.26 & 44.21 \\
2700 & 74.96 & 22.43 & 57.02 \\
3240 & 77.19 & 27.54 & 67.92 \\
3780 & 84.18 & 32.85 & 71.63 \\
4320 & 172.57 & 38.88 & 175.75 \\
4860 & 214.66 & 43.48 & 235.29 \\
5400 & 317.09 & 51.96 & 255.28 \\
5940 & 391.10 & 58.71 & 320.05 \\
6480 & 466.04 & 66.52 & 478.27 \\
7020 & 555.22 & 77.08 & 557.79 \\
7560 & 706.52 & 93.23 & 643.21 \\
8100 & 718.93 & 127.17 & 752.10 \\
8640 & 908.87 & 135.06 & 796.29 \\
9180 & 962.22 & 163.37 & 932.82 \\
9720 & 1007.20 & 195.49 & 1057.00 \\
10260 & 1247.08 & 217.08 & 1106.27 \\
10800 & 1571.24 & 255.97 & 1207.11 \\
11340 & 1853.91 & 290.33 & 1393.50 \\
11880 & 2095.26 & 327.06 & 1832.32 \\
12420 & 2240.76 & 382.05 & 1946.83 \\
12960 & 2417.49 & 412.56 & 2044.38 \\
\hline
\end{tabular}
\caption{Comparison of Total Runtimes for Adaptive Weights using \texttt{QuadCon}, 1-layer and 2-layer McCormick}
\label{tab:comparison_total_runtimes_aws}
\end{table}

The runtime for the \texttt{QuadCon} method initially increases polynomially, but at higher $n$ values, it begins to show signs of exponential growth, highlighting the computational challenges of handling quadratic constraints in higher dimensions. To assess solution quality, the Hypervolume (HV) Indicator is used, with comparisons for different methods presented in Table \ref{tab:hv_comparison_methods_aws}, using the same reference point for each $n$.
\begin{table}[h!]
    \centering
    \begin{tabular}{ccccc}
        \toprule
        \textbf{n} & \textbf{Pareto Points} & \textbf{HV (Q) ($\cdot 10^7$)} & \textbf{HV (1-l. M) ($\cdot 10^7$)} & \textbf{HV (2-l. M) ($\cdot 10^7$)} \\
        \midrule
        540 & 6 & 0.22 & 0.23 & 0.23 \\
        1080 & 6 & 0.88 & 0.91 & 0.91 \\
        1620 & 6 & 1.99 & 2.05 & 2.05 \\
        2160 & 6 & 3.51 & 3.65 & 3.65 \\
        2700 & 6 & 5.56 & 5.75 & 5.75 \\
        \bottomrule
    \end{tabular}
    \caption{Comparison of Hypervolume Indicator for Adaptive Weights using \texttt{QuadCon}, 1-layer McCormick, and 2-layer McCormick methods}
    \label{tab:hv_comparison_methods_aws}
\end{table}
The McCormick relaxation methods (single-layer and 2-layer) reduce computational burden without compromising the quality of Pareto solutions compared to \texttt{QuadCon}. The slight differences in Hypervolume (HV) values are due to the broader exploration of the Pareto front enabled by relaxing quadratic constraints, allowing for a larger search space. Overall, the McCormick methods offer a significant improvement in runtime while maintaining similar solution quality to \texttt{QuadCon}, with the single-layer McCormick method being the most efficient for the MOUC problem, balancing both efficiency and solution quality.

\subsection{Epsilon Constraints}
In this section, we explore the epsilon constraints method for solving our biobjective problem of interest. As earlier, we focus on both versions, i.e. \texttt{Quadcon} and McCormick relaxations.

\subsubsection{Quadratic Constraints with \texttt{QuadCon}}
First, we discuss the process of optimization using varying $\varepsilon$ for the quadratic constraints using Gurobi's \texttt{QuadCon} functionality. The corresponding formulation has been discussed in the previous section.
As earlier, we vary the problem dimension $n$ and record the performance of Gurobi and \texttt{Quadcon}. For $n=2160$, Table \ref{tab:optimization_results_f1} presents the results of optimizing the first objective function, \( f_1 \), while treating the second objective function \( f_2 \) as a constraint, with varying $\varepsilon$ values. The upper bound for \( f_2 \) increases incrementally with the values of $\varepsilon$ (from 0 to 1). 
\begin{table}[H]
\centering
\begin{tabular}{ccc}
\hline
\textbf{\(\varepsilon\)} & \textbf{Upper Bound for \(f_2\)} & \textbf{Optimal \(f_1(x)\)} \\
\hline
0.0 & \(f_2(x) \leq 12.544\) & 3080916.175 \\
0.1 & \(f_2(x) \leq 23.400\) & 2977178.832 \\
0.2 & \(f_2(x) \leq 34.255\) & 2940355.755 \\
0.3 & \(f_2(x) \leq 45.111\) & 2918969.545 \\
0.4 & \(f_2(x) \leq 55.966\) & 2898231.487 \\
0.5 & \(f_2(x) \leq 66.822\) & 2884766.524 \\
0.6 & \(f_2(x) \leq 77.677\) & 2871860.374 \\
0.7 & \(f_2(x) \leq 88.533\) & 2860818.942 \\
0.8 & \(f_2(x) \leq 99.388\) & 2851375.563 \\
0.9 & \(f_2(x) \leq 110.244\) & 2843023.920 \\
1.0 & \(f_2(x) \leq 121.099\) & 2835606.661 \\
\hline
\end{tabular}
\caption{Results for \(\varepsilon\)-constraints using \texttt{QuadCon}}
\label{tab:optimization_results_f1}
\end{table}
Tables \ref{tab:total_runtime_f1} and \ref{tab:gurobi_runtime_f1} provide an overview of the total runtime and the Gurobi runtime, respectively, for increasing $\varepsilon$ values from 0 to 1, and increasing dimensionality ($n$). 
From Tables \ref{tab:total_runtime_f1} and \ref{tab:gurobi_runtime_f1}, we observe that the total and Gurobi runtimes increase with the number of decision variables and the decrease in step sizes. For smaller $n$ (say, $n = 2160$), the difference between the total runtime and Gurobi runtime is relatively small. This indicates that a significant portion of the runtime is spent on actual optimization. However, as the number of decision variables increases (e.g., 10800), the proportion of total runtime to Gurobi runtime increases in proportion. This difference can be attributed to the overhead of data pre-processing and construction of objective and constraint matrices, which become substantial as the problem size grows.
\begin{table}[h!]
\centering
\begin{tabular}{|c|c|c|c|c|c|}
\hline
\diagbox{\textbf{n}}{\textbf{Step Size}} & \textbf{0.2} & \textbf{0.1} & \textbf{0.05} & \textbf{0.025} & \textbf{0.0125} \\ \hline
\textbf{1080} & 4.52 & 9.61 & 17.83 & 33.3 & 69.43 \\ \hline
\textbf{2160} & 7.41 & 15.77 & 30.85 & 56.66 & 111.94 \\ \hline
\textbf{3240} & 10.75 & 26.79 & 49.31 & 90.78 & 162.42 \\ \hline
\textbf{4320} & 35.42 & 63.84 & 116.8 & 180.44 & 307.76 \\ \hline
\textbf{5400} & 48.49 & 84.96 & 146.7 & 252.59 & 443.45 \\ \hline
\textbf{6480} & 75.04 & 122.37 & 180.1 & 290.77 & 545.1 \\ \hline
\textbf{7560} & 115.76 & 166.78 & 243.61 & 414.27 & 697.38 \\ \hline
\textbf{8640} & 145.07 & 240.78 & 349.66 & 464.79 & 888.99 \\ \hline
\textbf{9720} & 228.25 & 334.45 & 417.19 & 591.14 & 1088.37 \\ \hline
\textbf{10800} & 314.4 & 401.66 & 501.33 & 687.98 & 1228.2 \\ \hline
\end{tabular}
\caption{Total Runtime (s) for \(\varepsilon\)-constraints using \texttt{QuadCon} with objective \( f_1 \) and constraining $f_2$}
\label{tab:total_runtime_f1}
\end{table}
\begin{table}[h!]
\centering
\begin{tabular}{|c|c|c|c|c|c|}
\hline
\diagbox{\textbf{n}}{\textbf{Step Size}} & \textbf{0.2} & \textbf{0.1} & \textbf{0.05} & \textbf{0.025} & \textbf{0.0125} \\ \hline
\textbf{1080} & 3.71 & 8.73 & 16.92 & 32.08 & 67.57 \\ \hline
\textbf{2160} & 5.48 & 13.96 & 29.08 & 54.66 & 109.22 \\ \hline
\textbf{3240} & 6.44 & 22.21 & 44.4 & 85.12 & 155.06 \\ \hline
\textbf{4320} & 22.57 & 53.28 & 105.92 & 169.18 & 296.53 \\ \hline
\textbf{5400} & 23.07 & 61.43 & 122.93 & 230.42 & 418.48 \\ \hline
\textbf{6480} & 25.48 & 66.81 & 129.74 & 241.04 & 496.05 \\ \hline
\textbf{7560} & 26.77 & 90.03 & 170.84 & 335.35 & 612.15 \\ \hline
\textbf{8640} & 29.88 & 115.33 & 218.1 & 342.95 & 762.71 \\ \hline
\textbf{9720} & 57.36 & 140.43 & 241.4 & 422.32 & 925.16 \\ \hline
\textbf{10800} & 64.05 & 167.65 & 248.98 & 439.03 & 975.69 \\ \hline
\end{tabular}
\caption{Gurobi Runtime (s) for \(\varepsilon\)-constraints using \texttt{QuadCon} with objective \( f_1 \) and constraining $f_2$}
\label{tab:gurobi_runtime_f1}
\end{table}
Table \ref{tab:hv_epsilon_quadcon_f1} provides the Hypervolume (HV) Indicator values for $n = 2160$ and different step sizes in the $\varepsilon$-constraints method for \( f_1 \) using \texttt{QuadCon}. The number of Pareto points obtained for each step size is also indicated in table~\ref{tab:hv_epsilon_quadcon_f1}.
\begin{table}[h!]
    \centering
    \begin{tabular}{cccccc}
        \toprule
        \textbf{Step Size} & \textbf{Total (s)} & \textbf{Gurobi (s)} & \textbf{Pareto Points} & \textbf{HV ($\cdot 10^7$)} \\
        \midrule
        0.2 & 7.41 & 5.48 & 6 & 4.29 \\
        0.1 & 15.77 & 13.96 & 11 & 4.46 \\
        0.05 & 30.85 & 29.08 & 21 & 4.54 \\
        0.025 & 56.66 & 54.66 & 41 & 4.58 \\
        0.0125 & 111.94 & 109.22 & 81 & 4.60 \\
        \bottomrule
    \end{tabular}
    \caption{Hypervolume Indicator for \(\varepsilon\)-constraints using \texttt{QuadCon} with objective \( f_1 \) and constraining $f_2$}
    \label{tab:hv_epsilon_quadcon_f1}
\end{table}
The number of Pareto points and the Hypervolume (HV) values increase as step sizes decrease, due to finer granularity that allows for a more detailed exploration of the Pareto front. This leads to a better representation of the Pareto front, although it also requires more computational effort. Additionally, the study also examined the optimization of $f_2$ with a constraint on $f_1$. The total and Gurobi runtime for different step sizes and numbers of decision variables are shown in Tables \ref{tab:total_runtime_f2} and \ref{tab:gurobi_runtime_f2}, respectively.
\begin{table}[h!]
\centering
\begin{tabular}{|c|c|c|c|c|c|}
\hline
\diagbox{\textbf{n}}{\textbf{Step Size}} & \textbf{0.2} & \textbf{0.1} & \textbf{0.05} & \textbf{0.025} & \textbf{0.0125} \\ \hline
\textbf{1080} & 12.06 & 20.12 & 28.76 & 64.57 & 130.56 \\ \hline
\textbf{2160} & 19.22 & 33.36 & 57.72 & 141.72 & 268.95 \\ \hline
\textbf{3240} & 77.99 & 94.41 & 204.25 & 388.61 & 719.98 \\ \hline
\textbf{4320} & 141.28 & 169.26 & 414.99 & 753.04 & 1603.93 \\ \hline
\textbf{5400} & 315.77 & 375.54 & 956.30 & 1909.03 & 3691.47 \\ \hline
\textbf{6480} & 808.99 & 918.20 & 3386.30 & 5680.24 & 10752.95 \\ \hline
\textbf{7560} & 1501.76 & 1685.81 & 7416.45 & 15472.55 & * \\ \hline
\textbf{8640} & 3287.61 & 3531.25 & 18339.38 & * & * \\ \hline
\textbf{9720} & 12822.71 & 13182.87 & * & * & * \\ \hline
\textbf{10800} & * & * & * & * & * \\ \hline
\end{tabular}
\caption{Total Runtime (s) for \(\varepsilon\)-constraints using \texttt{QuadCon} with objective \( f_2 \) and constraining $f_1$}
\label{tab:total_runtime_f2}
\end{table}
\begin{table}[h!]
\centering
\begin{tabular}{|c|c|c|c|c|c|}
\hline
\diagbox{\textbf{n}}{\textbf{Step Size}} & \textbf{0.2} & \textbf{0.1} & \textbf{0.05} & \textbf{0.025} & \textbf{0.0125} \\ \hline
\textbf{1080} & 10.48 & 19.32 & 27.87 & 63.38 & 128.91 \\ \hline
\textbf{2160} & 17.33 & 30.98 & 55.30 & 139.10 & 266.12 \\ \hline
\textbf{3240} & 72.29 & 89.37 & 199.89 & 383.75 & 712.82 \\ \hline
\textbf{4320} & 130.52 & 159.04 & 403.29 & 741.02 & 1591.85 \\ \hline
\textbf{5400} & 290.62 & 352.39 & 933.15 & 1884.49 & 3666.34 \\ \hline
\textbf{6480} & 760.29 & 875.91 & 3344.42 & 5632.03 & 10693.24 \\ \hline
\textbf{7560} & 1423.18 & 1619.74 & 7349.72 & 15402.42 & * \\ \hline
\textbf{8640} & 3185.15 & 3425.45 & 18237.44 & * & * \\ \hline
\textbf{9720} & 12664.94 & 13042.98 & * & * & * \\ \hline
\textbf{10800} & * & * & * & * & * \\ \hline
\end{tabular}
\caption{Gurobi Runtime (s) for \(\varepsilon\)-constraints using \texttt{QuadCon} with objective \( f_2 \) and constraining $f_1$}
\label{tab:gurobi_runtime_f2}
\end{table}
The runtime trends for optimizing $f_2$ follow a similar pattern to those for $f_1$, with both total runtime and Gurobi runtime increasing as the number of decision variables and the step size decrease. However, the increase in runtime is more pronounced for $f_2$, especially for larger values of $n$. For example, the total runtime exceeds 12,822 seconds for $n = 9720$ and step size 0.2, and surpasses 20,000 seconds for smaller step sizes and larger $n$. Additionally, the proportion of Gurobi runtime within the total runtime is higher when optimizing $f_2$, indicating that a significant portion of the time is spent on solving the optimization problem using Gurobi. The longer runtime for $f_2$ is attributed to the increased computational complexity of constraining $f_1$, which results in more complex constraints and a larger objective space.
Table \ref{tab:hv_epsilon_quadcon_f2} presents the Hypervolume (HV) values for the $\varepsilon$-constraints method when optimizing \( f_2 \) using \texttt{QuadCon}, across various step sizes.
\begin{table}[h!]
    \centering
    \begin{tabular}{cccccc}
        \toprule
        \textbf{Step Size} & \textbf{Total (s)} & \textbf{Gurobi (s)} & \textbf{Pareto Points} & \textbf{HV ($\cdot 10^7$)} \\
        \midrule
        0.2 & 19.22 & 17.33 & 6 & 2.43 \\
        0.1 & 33.36 & 30.98 & 11 & 2.58 \\
        0.05 & 57.72 & 55.30 & 21 & 2.66 \\
        0.025 & 141.72 & 139.10 & 41 & 2.69 \\
        0.0125 & 268.95 & 266.12 & 81 & 2.71 \\
        \bottomrule
    \end{tabular}
    \caption{Hypervolume Indicator for \(\varepsilon\)-constraints using \texttt{QuadCon} with objective \( f_2 \) and constraining $f_1$}
    \label{tab:hv_epsilon_quadcon_f2}
\end{table}
Similar to the results for $f_1$, the number of distinct Pareto points increases with smaller step sizes, as finer granularity allows for more detailed exploration and a higher number of optimization runs. This leads to an increase in HV values and an enhanced Pareto frontier.

\subsubsection{Quadratic Constraints with McCormick}
With the introduction of McCormick relaxations, the number of decision variables is seen to increase from $n$ to $3n$. Note that this is due to the introduction of auxiliary variables $y_2$ and $w_2$.
For $n=2160$, optimization results for the objective $f_1$ using McCormick relaxations are presented in Table \ref{tab:optimization_results_f1_mccormick}. The table shows both the upper bounds and actual values of $f_2$ for varying $\varepsilon$ values
along with the optimal values of $f_1$. The results highlight that McCormick relaxation sometimes leads to deviations from the upper bound for $f_2$, which is expected. Adjustments to the bounds are made to ensure feasibility and to avoid generating identical solutions. 
As seen in the results with \texttt{QuadCon} in Table \ref{tab:optimization_results_f1}, as $\varepsilon$ increases from 0.0 to 1.0, the optimal value of $f_1$ decreases, suggesting that relaxing the constraint in $f_2$ allows for a more favorable optimization of $f_1$. Similarly, both the upper bound and actual value of $f_2$ increase as $\varepsilon$ increases, reflecting the relaxation of the constraint on $f_2$.
\begin{table}[H]
\centering
\begin{tabular}{|c|c|c|c|}
\hline
\(\varepsilon\) & Upper Bound for \(f_2(x)\) & \(f_2(x)\) & Optimal \(f_1(x)\) \\
\hline
0.0 & 12.544 & 13.147 & 3070051.332 \\
0.1 & 23.508 & 24.085 & 2974187.435 \\
0.2 & 34.472 & 35.069 & 2938771.839 \\
0.3 & 45.435 & 46.089 & 2917024.659 \\
0.4 & 56.399 & 57.143 & 2896820.000 \\
0.5 & 67.363 & 68.126 & 2883210.770 \\
0.6 & 78.327 & 79.197 & 2870070.570 \\
0.7 & 89.290 & 90.178 & 2859408.273 \\
0.8 & 100.254 & 101.182 & 2849989.210 \\
0.9 & 111.218 & 112.233 & 2841494.910 \\
1.0 & 122.182 & 122.182 & 2834788.460 \\
\hline
\end{tabular}
\caption{Results for objective $f_1$ and \(\varepsilon\)-constraints on $f_2$ using McCormick}
\label{tab:optimization_results_f1_mccormick}
\end{table}
Tables \ref{tab:epsilon_mcormick_total_f1} and \ref{tab:epsilon_mcormick_gurobi_f1} present the total and Gurobi runtime, respectively, for varying $n$ and step sizes using the single-layer McCormick relaxation. The number of decision variables presented only accounts for the original decision variables \( x \). 
\begin{table}[h!]
\centering
\begin{tabular}{|c|c|c|c|c|c|c|}
\hline
\diagbox{\textbf{n}}{\textbf{Step Size}} & \textbf{0.2} & \textbf{0.1} & \textbf{0.05} & \textbf{0.025} & \textbf{0.0125} \\ \hline
\textbf{1080} & 4.67 & 9.2 & 17.41 & 35.03 & 70.64 \\ \hline
\textbf{2160} & 9.01 & 14.75 & 25.13 & 45.57 & 85.11 \\ \hline
\textbf{3240} & 12.83 & 22.35 & 36.03 & 62.21 & 123.16 \\ \hline
\textbf{4320} & 19.95 & 35.5 & 53.05 & 83.83 & 142.29 \\ \hline
\textbf{5400} & 32.31 & 49.08 & 70.49 & 112.26 & 191.23 \\ \hline
\textbf{6480} & 47.2 & 61.64 & 99.77 & 147.93 & 247.26 \\ \hline
\textbf{7560} & 78.43 & 98.27 & 149.09 & 198.01 & 335.08 \\ \hline
\textbf{8640} & 92.67 & 123.13 & 165.25 & 259.01 & 389.19 \\ \hline
\textbf{9720} & 119.39 & 168.26 & 231.31 & 317.09 & 464.33 \\ \hline
\textbf{10800} & 144.59 & 193.57 & 254.91 & 370.33 & 589.4 \\ \hline
\end{tabular}
\caption{Total Runtime (s) for objective $f_1$ and \(\varepsilon\)-constraints on $f_2$ using McCormick}
\label{tab:epsilon_mcormick_total_f1}
\end{table}

\begin{table}[h!]
\centering
\begin{tabular}{|c|c|c|c|c|c|c|}
\hline
\diagbox{\textbf{n}}{\textbf{Step Size}} & \textbf{0.2} & \textbf{0.1} & \textbf{0.05} & \textbf{0.025} & \textbf{0.0125} \\ \hline
\textbf{1080} & 3.84 & 8.17 & 16.27 & 33.87 & 68.67 \\ \hline
\textbf{2160} & 5.67 & 12.16 & 22.54 & 42.79 & 81.99 \\ \hline
\textbf{3240} & 5.83 & 15.54 & 29.81 & 55.08 & 109.8 \\ \hline
\textbf{4320} & 6.02 & 22.03 & 39.58 & 69.43 & 128.36 \\ \hline
\textbf{5400} & 9.39 & 26.02 & 45.61 & 83.79 & 167.55 \\ \hline
\textbf{6480} & 11.12 & 27.47 & 65.61 & 113.47 & 210.19 \\ \hline
\textbf{7560} & 13.26 & 37.36 & 93.12 & 143.08 & 281.36 \\ \hline
\textbf{8640} & 16.88 & 47.98 & 93.76 & 185.1 & 316.61 \\ \hline
\textbf{9720} & 18.77 & 57.21 & 120.53 & 219.87 & 364.02 \\ \hline
\textbf{10800} & 14.24 & 62.49 & 130.73 & 236.99 & 458.01 \\ \hline
\end{tabular}
\caption{Gurobi Runtime (s) for objective $f_1$ and \(\varepsilon\)-constraints on $f_2$ using McCormick}
\label{tab:epsilon_mcormick_gurobi_f1}
\end{table}
As expected, we observe that both the total and Gurobi runtimes increase with the number of decision variables and the decrease in step sizes. Comparing these results with those obtained using \texttt{QuadCon}, we see that the McCormick relaxation method improves the runtime performance, particularly with larger problem dimension. The reduction in runtime can be attributed to the relaxation of quadratic constraints into linear constraints, which are computationally less intensive for the solver to handle, despite the addition of auxiliary variables.
Table \ref{tab:hv_epsilon_mccormick_f1} provides the Hypervolume (HV) Indicator values and the number of Pareto points for $n = 2160$ using different step sizes in the $\varepsilon$-constraints method for \( f_1 \) with \texttt{QuadCon}. As expected, while smaller step sizes require more computational effort, they yield more Pareto points and higher HV values, indicating an improved quality of the Pareto front.
\begin{table}[h!]
    \centering
    \begin{tabular}{cccccc}
        \toprule
        \textbf{Step Size} & \textbf{Total (s)} & \textbf{Gurobi (s)} & \textbf{Pareto Points} & \textbf{HV ($\cdot 10^7$)} \\
        \midrule
        0.2 & 9.01 & 5.67 & 6 & 4.29 \\
        0.1 & 14.75 & 12.16 & 11 & 4.46 \\
        0.05 & 25.13 & 22.54 & 21 & 4.54 \\
        0.025 & 45.57 & 42.79 & 41 & 4.58 \\
        0.0125 & 85.11 & 81.99 & 81 & 4.60 \\
        \bottomrule
    \end{tabular}
    \caption{Hypervolume Indicator for \(\varepsilon\)-constraints for minimizing $f_1$ and constraining $f_2$ using McCormick}
    \label{tab:hv_epsilon_mccormick_f1}
\end{table}
Similarly, tables \ref{tab:epsilon_mccormick_total_f2} and \ref{tab:epsilon_mccormick_gurobi_f2} provide an overview of the total runtime and the Gurobi runtime, respectively, for different step sizes and increasing numbers of decision variables for optimizing \( f_2 \) and constraining $f_1$ using the McCormick relaxation method.
\begin{table}[h!]
\centering
\begin{tabular}{|c|c|c|c|c|c|}
\hline
\diagbox{\textbf{n}}{\textbf{Step Size}} & \textbf{0.2} & \textbf{0.1} & \textbf{0.05} & \textbf{0.025} & \textbf{0.0125} \\ \hline
\textbf{1080} & 3.91 & 7.68 & 13.34 & 26.29 & 56.18 \\ \hline
\textbf{2160} & 7.5 & 13.68 & 18.86 & 35.98 & 73.51 \\ \hline
\textbf{3240} & 16.05 & 21.33 & 31.41 & 53.42 & 94.95 \\ \hline
\textbf{4320} & 20.87 & 30.03 & 40.85 & 64.88 & 113.12 \\ \hline
\textbf{5400} & 32.85 & 40.51 & 55.09 & 83.68 & 137.94 \\ \hline
\textbf{6480} & 44.55 & 54.14 & 67.45 & 110.31 & 181.16 \\ \hline
\textbf{7560} & 67.02 & 79.04 & 94.29 & 122.47 & 201.13 \\ \hline
\textbf{8640} & 89.63 & 112.57 & 126.41 & 154.58 & 257.86 \\ \hline
\textbf{9720} & 116.9 & 132.33 & 144.78 & 186.97 & 293.91 \\ \hline
\textbf{10800} & 151.68 & 163.44 & 187 & 219.19 & 330.87 \\ \hline
\end{tabular}
\caption{Total Runtime (s) for minimizing $f_2$ and constraining $f_1$ using McCormick}
\label{tab:epsilon_mccormick_total_f2}
\end{table}
\begin{table}[h!]
\centering
\begin{tabular}{|c|c|c|c|c|c|}
\hline
\diagbox{\textbf{n}}{\textbf{Step Size}} & \textbf{0.2} & \textbf{0.1} & \textbf{0.05} & \textbf{0.025} & \textbf{0.0125} \\ \hline
\textbf{1080} & 3.03 & 6.72 & 12.14 & 25.39 & 54.94 \\ \hline
\textbf{2160} & 4.51 & 10.14 & 15.97 & 32.77 & 69.7 \\ \hline
\textbf{3240} & 6.1 & 13.11 & 22.99 & 45.14 & 87.32 \\ \hline
\textbf{4320} & 6.64 & 14.96 & 26.23 & 51.07 & 98.01 \\ \hline
\textbf{5400} & 6.22 & 15.52 & 27.47 & 56.81 & 112.38 \\ \hline
\textbf{6480} & 7.07 & 16.27 & 29.72 & 60.76 & 140.48 \\ \hline
\textbf{7560} & 10.82 & 21.38 & 35.79 & 66.7 & 146.63 \\ \hline
\textbf{8640} & 10.39 & 25.56 & 40.35 & 76.18 & 183.45 \\ \hline
\textbf{9720} & 12.57 & 29.48 & 42.58 & 85.65 & 197.26 \\ \hline
\textbf{10800} & 12.74 & 30.92 & 50.67 & 89.88 & 204.72 \\ \hline
\end{tabular}
\caption{Gurobi Runtime (s) for minimizing $f_2$ and constraining $f_1$ using McCormick}
\label{tab:epsilon_mccormick_gurobi_f2}
\end{table}
The comparative statistics between \texttt{Quadcon} and McCormick relaxations are similar here as well, where a significant reduction in the computational effort with the latter. For instance, with \( n = 10800 \) and a step size of 0.2, Gurobi solved the problem in 12.74 seconds using the McCormick relaxation, whereas with \texttt{QuadCon}, it exceeded $20,000$ seconds. Additionally, the total runtime with McCormick, including overhead, is 151.68 seconds, demonstrating the efficiency of the McCormick approach despite the additional auxiliary variables.
\begin{table}[h!]
    \centering
    \begin{tabular}{ccccc}
        \toprule
        \textbf{Step Size} & \textbf{Total (s)} & \textbf{Gurobi (s)} & \textbf{Pareto Points} & \textbf{HV ($\cdot 10^7$)} \\
        \midrule
        0.2 & 7.50 & 4.51 & 6 & 2.43 \\
        0.1 & 13.68 & 10.14 & 11 & 2.55 \\
        0.05 & 18.86 & 15.97 & 21 & 2.59 \\
        0.025 & 35.98 & 32.77 & 41 & 2.62 \\
        0.0125 & 73.51 & 69.70 & 81 & 2.63 \\
        \bottomrule
    \end{tabular}
    \caption{Hypervolume Indicator for \(\varepsilon\)-constraints for mimizing $f_2$ and constraining $f_1$ using McCormick}
    \label{tab:runtime_comparison_f2}
\end{table}

\subsubsection{\texttt{QuadCon} vs. McCormick}
\begin{figure}[h!]
\centering
\resizebox{0.65\textwidth}{!}{ 
\begin{tikzpicture}
    \begin{axis}[
        ybar,
        width=\textwidth,
        height=0.5\textheight,
        bar width=5pt,
        xlabel={n},
        ylabel={Runtime (s)},
        symbolic x coords={1080, 2160, 3240, 4320, 5400, 6480, 7560, 8640, 9720, 10800},
        xtick={1080, 2160, 3240, 4320, 5400, 6480, 7560, 8640, 9720, 10800},
        x tick label style={rotate=45, anchor=east}, 
        legend pos=north west,
        ymin=0,
        enlarge x limits=0.1
    ]
    \addplot coordinates {(1080,9.61) (2160,15.77) (3240,26.79) (4320,63.84) (5400,84.96) (6480,122.37) (7560,166.78) (8640,240.78) (9720,334.45) (10800,401.66)};
    \addplot coordinates {(1080,9.20) (2160,14.75) (3240,22.35) (4320,35.50) (5400,49.08) (6480,61.64) (7560,98.27) (8640,123.13) (9720,168.26) (10800,193.57)};
    \legend{\texttt{QuadCon} Runtime, McCormick Runtime}
    \end{axis}
\end{tikzpicture}
}
\caption{Comparison of Total Runtimes for \(\varepsilon\)-constraints using \texttt{QuadCon} vs. McCormick with Step Size 0.1 for optimizing $f_1$ and constraining $f_2$}
\label{fig:comparison_econst_total_f1}
\end{figure}
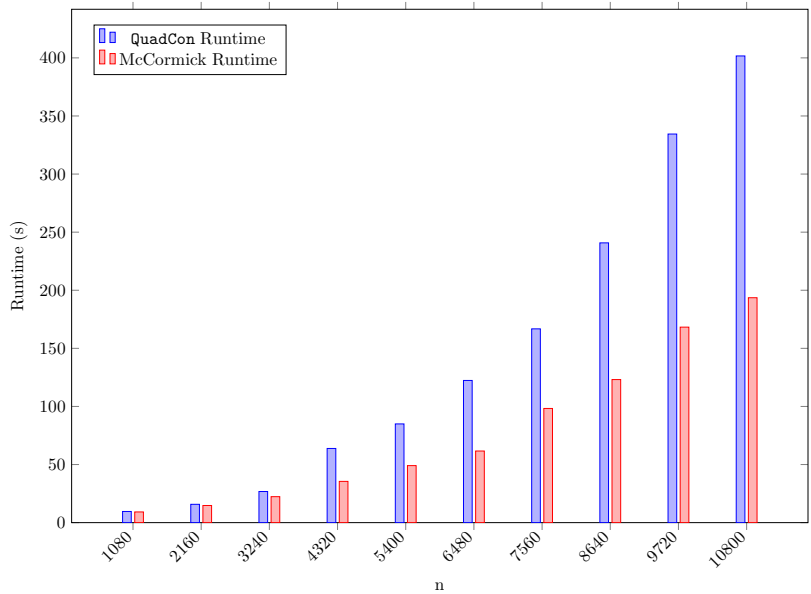

Figure \ref{fig:comparison_econst_total_f1} displays the total runtimes for the \texttt{QuadCon} and single-layer McCormick methods across various values of $n$. Note that in this case, $f_1$ is minimized and $f_2$ is constrained. These results further confirm that the single-layer McCormick method consistently demonstrates lower total runtimes, making it the faster approach for handling quadratic constraints compared to \texttt{QuadCon} in our optimization problem.
Table \ref{tab:hv_comparison_epilon_f1} compares the performance of the \texttt{QuadCon} and McCormick methods in the $\varepsilon$-constraints method for minimizing \( f_1 \) (and constraining $f_2$), presenting the number of Pareto points, Gurobi runtimes, and HV values calculated by using the same reference points for different step sizes. Despite the differences in runtimes, the McCormick relaxation method does not compromise on the quality of the Pareto solutions. This is reflected in the HV values, which remain consistently high and closely match those obtained with \texttt{QuadCon}.
\begin{table}[h!]
    \centering
    \begin{tabular}{ccccccc}
        \toprule
        \textbf{Step Size} & \textbf{Pareto Points} & \textbf{Runtime (Q)} & \textbf{Runtime (M)} & \textbf{HV (Q)} & \textbf{HV (M)} \\
        \midrule
        0.2 & 6 & 5.48 & 5.67 & 4.29 & 4.29 \\
        0.1 & 11 & 13.96 & 12.16 & 4.46 & 4.45 \\
        0.05 & 21 & 29.08 & 22.54 & 4.54 & 4.53 \\
        0.025 & 41 & 54.66 & 42.79 & 4.58 & 4.57 \\
        0.0125 & 81 & 109.22 & 81.99 & 4.60 & 4.59 \\
        \bottomrule
    \end{tabular}
    \caption{Comparison of Hypervolume Indicator for \(\varepsilon\)-constraints using \texttt{QuadCon} vs. McCormick for minimizing $f_1$ and contraining $f_2$}
    \label{tab:hv_comparison_epilon_f1}
\end{table}
Table \ref{tab:total_runtime_comparison_f2} compares the total runtimes for solving the optimization problem with objective $f_2$ (and constraining $f_1$) using the $\varepsilon$-constraints method with \texttt{QuadCon} and McCormick relaxations. The analysis shows that, similar to $f_1$, McCormick relaxation outperforms \texttt{QuadCon}, with differences becoming even more significant for $f_2$. Specifically, the complexity of the \texttt{QuadCon} method grows exponentially as the number of decision variables increases. For larger values of $n$ (greater than 9720), \texttt{QuadCon} runtimes exceed 20,000 seconds, while McCormick relaxation allows the problem to be solved in just 163.44 seconds.
\begin{table}[h!]
\centering
\begin{tabular}{|c|c|c|}
\hline
\textbf{n} & \textbf{\texttt{QuadCon} Runtime(s)} & \textbf{McCormick Runtime(s)} \\ \hline
\textbf{1080} & 20.12 & 7.68 \\ \hline
\textbf{2160} & 33.36 & 13.68 \\ \hline
\textbf{3240} & 94.41 & 21.33 \\ \hline
\textbf{4320} & 169.26 & 30.03 \\ \hline
\textbf{5400} & 375.54 & 40.51 \\ \hline
\textbf{6480} & 918.20 & 54.14 \\ \hline
\textbf{7560} & 1685.81 & 79.04 \\ \hline
\textbf{8640} & 3531.25 & 112.57 \\ \hline
\textbf{9720} & 13182.87 & 132.33 \\ \hline
\textbf{10800} & * & 163.44 \\ \hline
\end{tabular}
\caption{Comparison of Total Runtimes for \(\varepsilon\)-constraints using \texttt{QuadCon} vs. McCormick with Step Size 0.1 for minimizing $f_2$ and constraining $f_1$}
\label{tab:total_runtime_comparison_f2}
\end{table}
Table \ref{tab:combined_runtime_comparison_f2} compares the efficiency of \texttt{QuadCon} and McCormick in the $\varepsilon$-constraints method for minimizing $f_2$ (and constraining $f_1$), presenting the number of Pareto points, Gurobi runtimes, and HV values for different step sizes using the same reference points. Although the runtimes differ significantly between the two methods, the HV values are similar. The difference in HV is due to the constraint on $f_1$, when optimizing $f_2$, with the relaxation of $f_1$'s constraints resulting in a slightly different distribution of Pareto solutions, particularly focusing on regions with higher $f_1$ and lower $f_2$. Using tighter McCormick envelopes could potentially improve the HV, offering a more accurate representation of the Pareto front.
\begin{table}[h!]
    \centering
    \begin{tabular}{ccccccc}
        \toprule
        \textbf{Step Size} & \textbf{Pareto Points} & \textbf{Runtime (Q)} & \textbf{Runtime (M)} & \textbf{HV (Q)} & \textbf{HV (M)} \\
        \midrule
        0.2 & 6 & 17.33 & 4.51 & 2.43 & 2.43 \\
        0.1 & 11 & 30.98 & 10.14 & 2.58 & 2.55 \\
        0.05 & 21 & 55.30 & 15.97 & 2.66 & 2.59 \\
        0.025 & 41 & 139.10 & 32.77 & 2.69 & 2.62 \\
        0.0125 & 81 & 266.12 & 69.70 & 2.71 & 2.63 \\
        \bottomrule
    \end{tabular}
    \caption{Comparison of Hypervolume Indicator for \(\varepsilon\)-constraints using \texttt{QuadCon} vs. McCormick for minimizing $f_2$ and constraining $f_1$}
    \label{tab:combined_runtime_comparison_f2}
\end{table}
\subsection{Visualization of Pareto Frontiers:}
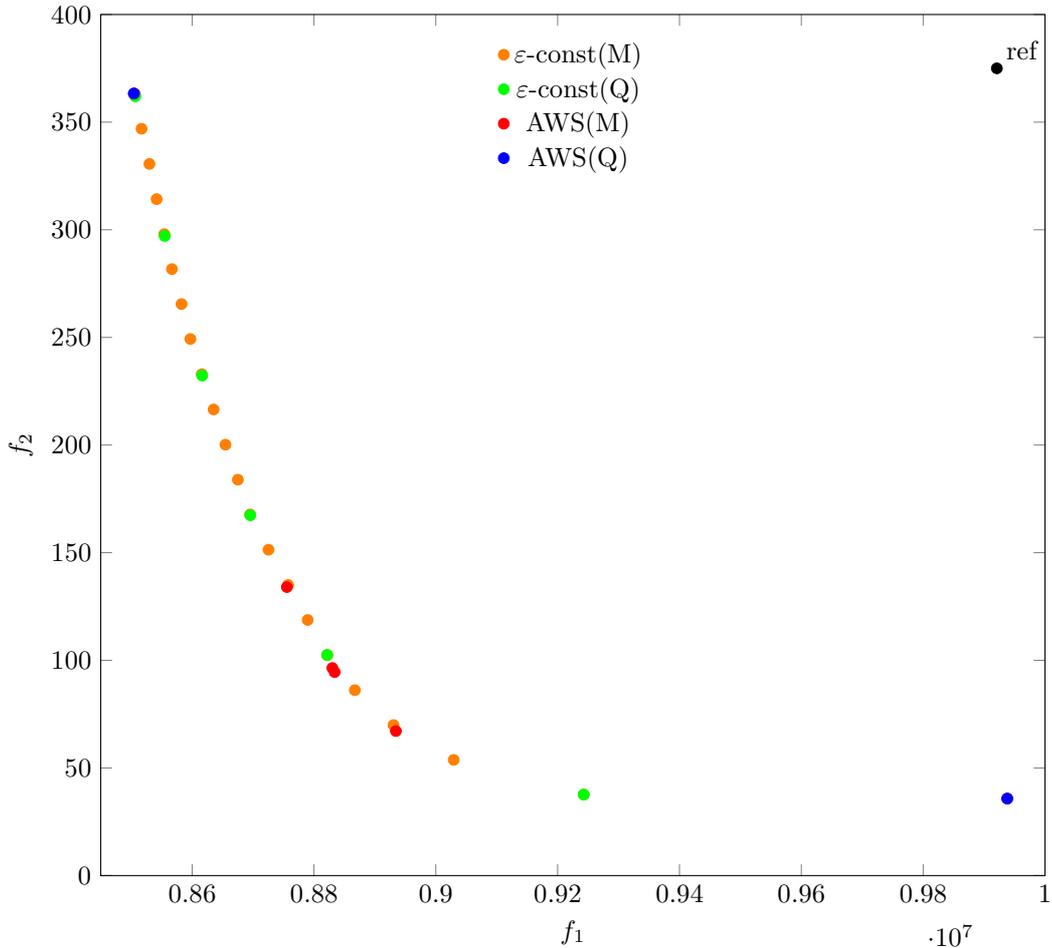
\begin{figure}[h]
    \centering
    \begin{tikzpicture}
        \begin{axis}[
            width=14cm,
            height=13cm,
            xlabel={$f_1$},
            ylabel={$f_2$},
            xmin=8450000, xmax=10000000,
            ymin=0, ymax=400,
            legend pos=outer north east,
            legend style={draw=none, at={(0.5,0.98)}, anchor=north},
            scatter/classes={d={mark=*,orange}, c={mark=*,green}, b={mark=*,red}, a={mark=*,blue}, e={mark=*,black}
            }
        ]
        \addplot[scatter, only marks, scatter src=explicit symbolic]
        table[meta=label] {
            x y label
            9242745.699 37.6315 d
            9029392.74 53.75 d
            8930490.148 69.9207 d
            8867095.968 86.1419 d
            8821863.482 102.3969 d
            8789703.62 118.7404 d
            8757543.759 135.0121 d
            8725383.897 151.3664 d
            8695146.475 167.7429 d
            8675074.066 183.9125 d
            8654655.496 200.166 d
            8635223.69 216.531 d
            8615791.883 232.893 d
            8597093.34 249.2562 d
            8582533.002 265.512 d
            8566859.97 281.7292 d
            8554295.886 297.9461 d
            8541731.802 314.245 d
            8529933.297 330.5785 d
            8516986.468 346.9351 d
            8506016.172 363.2276 d
        };
        \addlegendentry{$\varepsilon$-const(M)}
        \addplot[scatter, only marks, scatter src=explicit symbolic]
        table[meta=label] {
            x y label
            9242751.527 37.632 c
            8821584.214 102.5056 c
            8695519.761 167.3792 c
            8616503.417 232.2528 c
            8554928.999 297.1264 c
            8506720.42 362 c
        };
        \addlegendentry{$\varepsilon$-const(Q)}
        \addplot[scatter, only marks, scatter src=explicit symbolic]
        table[meta=label] {
            x y label
            9938102.025 35.76427583 b
            8504365.38 363.2967549 b
            8934472.346 67.13030114 b
            8755609.04 134.0240195 b
            8833937.148 94.53050021 b
            8830169.885 96.4297498 b
        };
        \addlegendentry{AWS(M)}
        \addplot[scatter, only marks, scatter src=explicit symbolic]
        table[meta=label] {
            x y label
            9938102.025 35.76427583 a
            8504365.38 363.2967549 a
        };
        \addlegendentry{AWS(Q)}
        \addplot[scatter, only marks, scatter src=explicit symbolic]
        table[meta=label] {
            x y label
            9921000 375 e
        };
        \node[anchor=south west] at (axis cs:9921000,375) {ref};
        \end{axis}
    \end{tikzpicture}
    \caption{Pareto Points for Different Methods for \( n=6480 \) with a time cap of 100 seconds}
    \label{fig:pareto_points_methods_1}
\end{figure}

\begin{figure}[h]
    \centering
    \begin{tikzpicture}
        \begin{axis}[
            width=14cm,
            height=13cm,
            xlabel={$f_1$},
            ylabel={$f_2$},
            xmin=14100000, xmax=16700000,
            ymin=0, ymax=660,
            legend pos=outer north east,
            legend style={draw=none, at={(0.5,0.98)}, anchor=north},
            scatter/classes={d={mark=*,orange}, c={mark=*,green}, b={mark=*,red}, a={mark=*,blue}, e={mark=*,black}
            }
        ]
        \addplot[scatter, only marks, scatter src=explicit symbolic]
        table[meta=label] {
            x y label
            15404576.06 62.7192 d
            15048987.7 89.5837 d
            14883888.84 116.5339 d
            14777779.42 143.5698 d
            14703105.75 170.6382 d
            14649505.98 197.7875 d
            14596065.05 224.9364 d
            14542306.44 252.2905 d
            14491910.76 279.5716 d
            14459315.73 306.5006 d
            14424425.79 333.61 d
            14392039.45 360.8052 d
            14359653.11 388.1547 d
            14328488.87 415.4271 d
            14304221.65 442.5189 d
            14278099.93 469.5487 d
            14257159.79 496.5985 d
            14236219.65 523.7413 d
            14215279.51 551.0322 d
            14194977.43 578.2251 d
            14176693.61 605.3793 d
        };
        \addlegendentry{$\varepsilon$-const(M)}
        \addplot[scatter, only marks, scatter src=explicit symbolic]
        table[meta=label] {
            x y label
            15404576.66 62.72 c
            14423837.87 334.1071 c
            14176631.14 605.4949 c
        };
        \addlegendentry{$\varepsilon$-const(Q)}
        \addplot[scatter, only marks, scatter src=explicit symbolic]
        table[meta=label] {
            x y label
            16570032.54 59.60712638 b
            14173942.3 605.4945914 b
            14892755.35 111.2781714 b
            14592666.09 223.3733659 b
            14724091.17 157.1159645 b
            14717624.98 160.3758666 b
        };
        \addlegendentry{AWS(M)}
        \addplot[scatter, only marks, scatter src=explicit symbolic]
        table[meta=label] {
            x y label
            16570032.54 59.60712638 a
            14173942.3 605.4945914 a
        };
        \addlegendentry{AWS(Q)}
        \addplot[scatter, only marks, scatter src=explicit symbolic]
        table[meta=label] {
            x y label
            16535000 625 e
        };
        \node[anchor=south west] at (axis cs:16535000,625) {ref};
        \end{axis}
    \end{tikzpicture}
    \caption{Pareto Points for Different Methods for \( n=10800 \) with a time cap of 300 seconds}
    \label{fig:pareto_points_methods_2}
\end{figure}
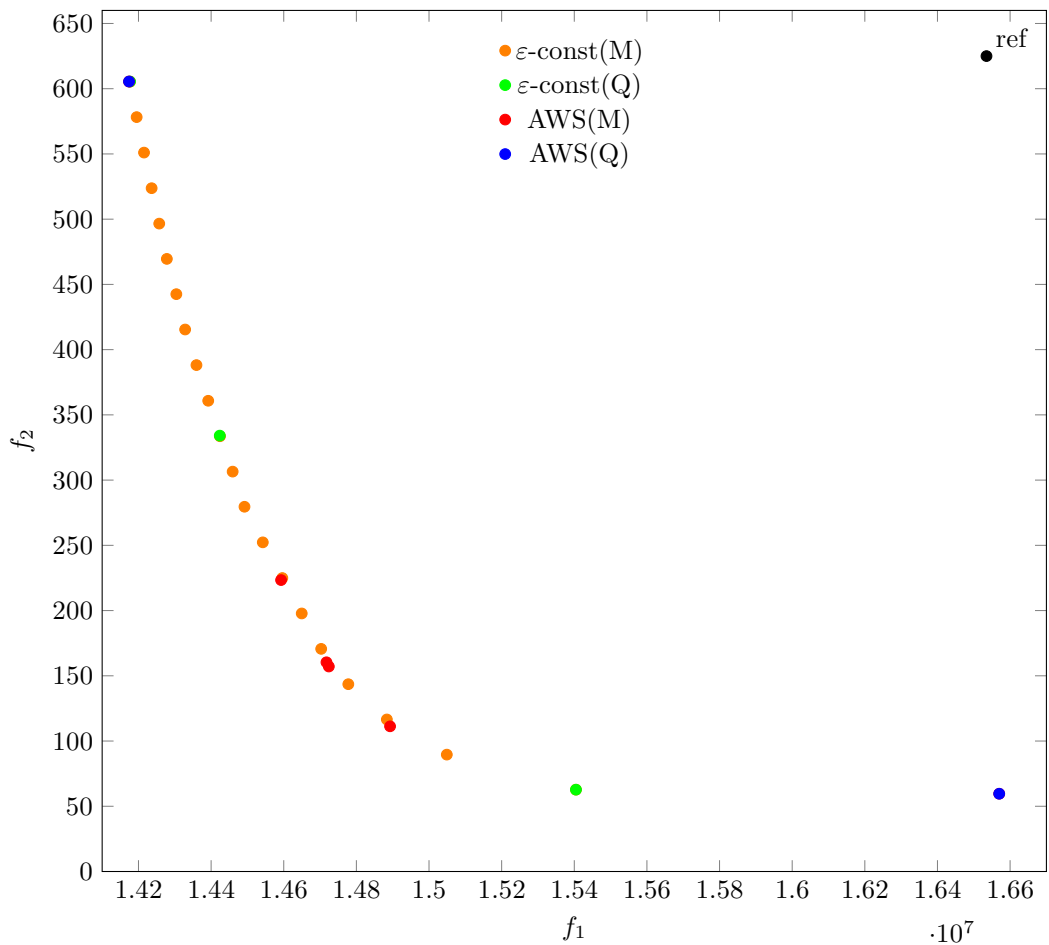

Figure \ref{fig:pareto_points_methods_1} compares the non-dominated solutions obtained using different methods: Adaptive Weighted Sum (AWS) with \texttt{QuadCon}, AWS with one-layer McCormick, $\varepsilon$-constraints with \texttt{QuadCon}, and $\varepsilon$-constraints with one-layer McCormick for $n=6480$, with a 100-second time cap on total runtime. The runtime includes overhead to prevent bias toward McCormick relaxation. Under this time cap, the results show: AWS with \texttt{QuadCon} yields 2 Pareto points, AWS with McCormick produces 6 points, $\varepsilon$-constraints with \texttt{QuadCon} gives 6 points, and $\varepsilon$-constraints with McCormick achieves 21 points. The Pareto solutions from each method are displayed, with most methods including the two extreme points of the Pareto front. Some points may overlap, making them less visible. Within each AWS and $\varepsilon$-constraints method, the \texttt{QuadCon} solutions are found in the McCormick results, but McCormick provides additional points within the same time limit.
Figure \ref{fig:pareto_points_methods_2} compares the Pareto points obtained for $n=10800$ with a 300-second total runtime cap, using the same methods: Adaptive Weighted Sum (AWS) with \texttt{QuadCon}, AWS with one-layer McCormick, $\varepsilon$-constraints with \texttt{QuadCon}, and $\varepsilon$-constraints with one-layer McCormick. Under this setup, AWS with \texttt{QuadCon} yields 2 Pareto points, AWS with McCormick provides 6 points, $\varepsilon$-constraints with \texttt{QuadCon} results in 3 points, and $\varepsilon$-constraints with McCormick generates 21 points. The solutions shown in the figure highlight the performance of each method. As in the previous analysis, the two extreme points are included in most methods, and within each AWS and $\varepsilon$-constraints approach, the \texttt{QuadCon} solutions are found in the McCormick results, which offer additional points. Overall, Figures \ref{fig:pareto_points_methods_1} and \ref{fig:pareto_points_methods_2} demonstrate that McCormick relaxation facilitates a more thorough exploration of the Pareto front, leading to a higher number of non-dominated solutions for a given computational budget.

\section{Conclusion}
\label{sec:conclusion}
In conclusion, this work has focused on improving the efficiency of solving the Multiobjective Unit Commitment (MOUC) problem, when tackled through integer programming approaches. We investigated the use of adaptive weighted methods and epsilon constraints for calculating Pareto frontiers, highlighting their advantages over uniform weights, despite the added complexity of quadratic constraints. To address this, we introduced McCormick relaxations to approximate the problem with linear constraints, leading to significant improvements in computational performance on a real-world case study. The results showed only a minimal trade-off in solution quality compared to the standard Gurobi solver applied to the original formulation. Based on these promising results, we propose further theoretical analysis of the McCormick-based algorithm and its potential extension to other application domains beyond unit commitment. While this study did not compare with genetic algorithms, future work could explore integrating mathematical programming with evolutionary methods to further enhance the solution process for such problems.
\section*{Acknowledgment}%
The research was funded partly by the Deutsche Forschungsgemeinschaft (DFG, German Research
Foundation) under Germany's Excellence Strategy --- The Berlin Mathematics
Research Center MATH+ (EXC-2046/1, project ID:\@390685689). The authors would additionally like to talk Prof. Falk Hante (Humboldt Universitaet zu Berlin) for insightful comments that helped with the technical aspects of this work.
\small
\bibliography{ref,sample,literature-ece}
\bibliographystyle{siam}
\end{document}

%% file: macros.tex
\usepackage{amsmath, amssymb, bbm, xspace}

\def\sign{\mathop{\hbox{\rm sign}}}

\def\spose#1{\hbox to 0pt{#1\hss}}

\def\text #1{\hbox{\quad#1\quad}}

\def\nthinsp{\mskip -2   mu}

\def\I{_{\scriptscriptstyle I}}

\def\Q{_{\scriptscriptstyle Q}}

\def\superstar{^{\raise 0.5pt\hbox{$\nthinsp *$}}}
\def\SUPERSTAR{^{\raise 0.5pt\hbox{$*$}}}

\def\lamstarT {\lambda^{\raise 0.5pt\hbox{$\nthinsp *$}T}}

\def\hbar{\skew{4.2}\bar h}

		\def\bkE{{\rm I\kern-.17em E}}
		\def\bk1{{\rm 1\kern-.17em l}}
		\def\bkD{{\rm I\kern-.17em D}}
		\def\bkR{{\rm I\kern-.17em R}}
		\def\bkP{{\rm I\kern-.17em P}}
		\def\bkY{{\bf \kern-.17em Y}}
		\def\bkZ{{\bf \kern-.17em Z}}

		\def\beq{\begin{eqnarray}}
		\def\bc{\begin{center}}
		\def\be{\begin{enumerate}}
		\def\bi{\begin{itemize}}
		\def\bs{\begin{small}}
		\def\bS{\begin{slide}}
		\def\ec{\end{center}}
		\def\ee{\end{enumerate}}
		\def\ei{\end{itemize}}
		\def\es{\end{small}}
		\def\eS{\end{slide}}
		\def\eeq{\end{eqnarray}}

	\def\cp2problem#1#2#3#4{\fbox
		 {\begin{tabular*}{0.9\textwidth}
			{@{}l@{\extracolsep{\fill}}l@{\extracolsep{6pt}}l@{\extracolsep{\fill}}c@{}}
				#1 & & $#4 $ 
			\end{tabular*}}}

\newcommand{\pmat}[1]{\begin{pmatrix} #1 \end{pmatrix}}
		
		\renewcommand{\emph}[1]{\textbf{#1}}

		\def\bkE{{\rm I\kern-.17em E}}
		\def\bk1{{\rm 1\kern-.17em l}}
		\def\bkD{{\rm I\kern-.17em D}}
		\def\bkR{{\rm I\kern-.17em R}}
		\def\bkP{{\rm I\kern-.17em P}}
		
		\def\bkZ{{\bf{Z}}}

\newcommand {\beeq}[1]{\begin{equation}\label{#1}}
\newcommand {\eeeq}{\end{equation}}
\newcommand {\bea}{\begin{eqnarray}}
\newcommand {\eea}{\end{eqnarray}}

\def\texitem#1{\par\smallskip\noindent\hangindent 25pt
               \hbox to 25pt {\hss #1 ~}\ignorespaces}

